\documentclass[reprint,superscriptaddress,aip,cha]{revtex4-1}
\usepackage{epsfig,graphicx,color}
\usepackage{amsmath, amsthm, amsfonts,amssymb}
\usepackage{dcolumn}
\usepackage{empheq}
\usepackage{tikz}
\makeatletter
\newcommand*{\centerfloat}{%
  \parindent \z@
  \leftskip \z@ \@plus 1fil \@minus \textwidth
  \rightskip\leftskip
  \parfillskip \z@skip}
\makeatother

\begin{document}

\title{Phase and amplitude dynamics of coupled oscillator systems on complex networks}\thanks{This article may be downloaded for personal use only. Any other use requires prior permission of the author and AIP Publishing. This article appeared in Chaos 30, 121102 (2020) and may be found at https://doi.org/10.1063/5.0031031.}

\author{Jae Hyung Woo}
\affiliation{Department of Psychological and Brain Sciences, Johns Hopkins University, Baltimore, MD 21218, USA}
\author{Christopher J. Honey}
\affiliation{Department of Psychological and Brain Sciences, Johns Hopkins University, Baltimore, MD 21218, USA}
\author{Joon-Young Moon}
 \email{joon.young.moon@gmail.com}
\affiliation{Department of Psychological and Brain Sciences, Johns Hopkins University, Baltimore, MD 21218, USA}

\date{\today}

\begin{abstract}
We investigated the locking behaviors of coupled limit-cycle oscillators with phase and amplitude dynamics. We focused on how the dynamics are affected by inhomogeneous coupling strength and by angular and radial shifts in the coupling function. We performed mean-field analyses of oscillator systems with inhomogeneous coupling strength, testing Gaussian, power-law, and brain-like degree  distributions. Even for oscillators with identical intrinsic frequencies and intrinsic amplitudes, we found that the coupling strength distribution and coupling function generated a wide repertoire of phase and amplitude dynamics. These included fully and partially locked states in which high-degree or low-degree nodes would phase-lead the network. The mean-field analytical findings were confirmed via numerical simulations. The results suggest that, in oscillator systems in which individual nodes can independently vary their amplitude over time, qualitatively different dynamics can be produced via shifts in the coupling strength distribution and the coupling form. Of particular relevance to information flows in oscillator networks, changes in the non-specific drive to individual nodes can make high-degree nodes phase-lag or phase-lead the rest of the network.

\end{abstract}
\maketitle

\begin{quotation}
Models of coupled oscillators have been widely used across variety of disciplines including physics, chemistry, and biology to describe the dynamics of systems with interacting elements. For example, fireflies adjust their blinking according to the light-flashes they see around them, and clusters of neurons produce rhythmic firing, whose timing depends on the input from other neurons. Previous studies have analyzed these coupled systems in phase-reduced models, in which each node in the network is treated like a cycling clock. However, in many real-world applications, such as in brain networks, the amplitude of activity at one location affects  the response at another location, just as more neuronal firing will produce a larger effect at recipient sites. Therefore, it is important to study dynamics of such systems using models that account for both the phase and amplitude of the dynamics of each node. To this end, we analyzed the dynamics of coupled identical oscillator  with phase and amplitude dynamics. We focused on how the dynamics are altered by two factors: first, the inhomogeneity in the coupling strength (so that some nodes have stronger connections than others) and second, the coupling function (how the response of a target node depends on the phase and amplitude of a source node). Building on previous works which have focused on the phase dimension of the dynamics, we mapped a rich repertoire of amplitude dynamics and phase dynamics, depending on the distribution of connection strengths and the form of the coupling function. Conditions for each of the possible classes of dynamics were identified using stability analysis, following a self-consistency argument, and the results were confirmed via numerical simulation. Of particular interest for neuroscience, we found that groups of nodes could shift from phase-leading to phase-lagging, depending on minor changes in the coupling function. Thus, small changes in a non-specific driving signal in the brain can cause shifts in the direction of signaling between brain regions. More generally, we also observed a variety of non-locked states ("drifting" dynamics) in which only a subset of strongly connected nodes in the network are strongly synchronized, while other nodes on the periphery operate more independently.

\end{quotation}

\section{\label{sec:intro}Introduction}
Coupled oscillator systems provide models for systems of interacting elements in many fields, including physics, chemistry, and biology.\cite{michaels1987mechanisms, pantaleone1998stability, ermentrout2001traveling, kuramoto2003chemical, ko2007effects, filatrella2008analysis, arenas2008synchronization, ferrari2015phase} The intrinsic dynamics of the oscillators, the couplings between the oscillators, and the connectivity among the oscillators jointly determine the dynamics of the coupled system. For the generality of the model and ease of analysis, many studies have focused on phase dynamics of the coupled oscillator systems.\cite{sakaguchi1987local, bonilla1992nonlinear, daido1992quasientrainment, nakamura1994clustering, kim1997multistability, strogatz2000kuramoto, chandra2019continuous} Even within phase-reduced systems, lacking any variation in amplitude, oscillator systems exhibit a rich repertoire of various synchronous behaviors such as in-phase synchronization, full locking, chimera state, and partial locking.\cite{chaos2019} However, the physics of many real-world systems includes dramatic amplitude dynamics: for example, in populations of neurons, the amplitude of oscillations varies dramatically across different brain regions and over time, with important functional implications. \cite{siegel_spectral_2012, de_pesters_alpha_2016, pfurtscheller_event-related_1999, jensen_oscillatory_2015} Therefore, oscillator models incorporating both phase and amplitude dynamics have been investigated, with a focus on the global stability of the system. These stability-focused studies revealed diverse outcomes, including amplitude death, chimera states, and phase-delay effects.\cite{ermentrout1990oscillator, matthews1990phase, choe2010controlling, d2010amplitude, laing2010chimeras, gambuzza2016amplitude, banerjee2018networks}. In the present work, rather than focusing on whether synchronization occurs, we focus on the variability in phase and amplitude across nodes, and how these are controlled by system-wide parameters. In particular, we focus on cases in which the network has inhomegeneous degree distributions, so that some nodes are more influential than others, and aim to characterize how high (or low) degree nodes may come to have increased (or decreased) amplitudes and to phase-lag (or phase-lead) the rest of the network.

The detailed interactions of phase-and-amplitude dynamics are especially important in the understanding of neural systems, because they may explain global state transitions associated with distinct modes of brain function. For example, Moon et al. \cite{moon2015general} showed that the phase-and-amplitude relationships differ before and after anesthesia, and demonstrated how this phenomenon could be captured in a Stuart-Landau model. More generally, it is critical to provide a theoretical framework to organize the diverse various locking behaviors observed in neural oscillatory systems, as the changes in these systems are associated with large-scale functional changes (sleep/wake, active/passive), and the controllers of such state-shifts remain largely unknown.\cite{honey2007network, cabral2011role, cabral2014exploring, sanz2015mathematical, finger2016modeling, breakspear2017dynamic, jymoon2017, vlasov2017hub, kim2018mechanisms, bansal2019cognitive, de2019connection, lee2019relationship, kang2020remote, kim2020alpha, woldman2020dynamic} \newline

In this paper, we study how the phase and amplitude dynamics of coupled oscillators depend on (i) spatial inhomogeneity in coupling strengths and (ii) the form of coupling function between nodes. These two factors combine to produce rich repertoire of synchronous behaviors in terms of both phases and amplitudes of the oscillators. Numerical simulations with Gaussian, power-law, and brain network distributions of the coupling inhomogeneity show the conditions for different categories of various synchronous states ranging from un-locked state, partially locked state to fully locked state. Analytical approximations were in agreement with our numerical results. This work can understood as a generalization of previous work examining coupling inhomogeneity in phase-based oscillators \cite{chaos2019}, extending the analysis to include variability in amplitude.

\section{\label{sec:derive}Model and Analysis}
We investigate a generalized form of coupled identical limit cycle oscillators, in the form of Stuart-Landau model. Stuart-Landau model is a canonical model, in the sense that oscillator systems reduce to the Stuart-Landau model near a Hopf bifurcation. \cite{guckenheimer1983nonlinear, hoppensteadt1997weakly}
Thus, we consider a mean-field model of $N$ weakly coupled Stuart-Landau oscillators, each near a Hopf bifurcation:
\begin{eqnarray} \label{eq:SL1}
\dot{z}_{j} &=& \{\lambda_j - |z_{j}|^2 + i\omega_j \}z_{j} + \frac{S K_{j}}{N}\sum_{k=1}^{N}(z_{k}e^{-i\beta}- z_{j}d_{0}e^{-i\alpha}), \nonumber \\
&&  j=1,2,...,N, ~\alpha \in [0,\pi) ,~\beta \in [0,\pi/2), ~d_0 \in \mathbb{R}, \nonumber \\
\end{eqnarray}
where $z_j(t)=r_j(t) e^{i\theta_j(t)}$ is the position of an oscillator $j$ in the complex plane at time $t$. $S$ is a parameter controlling the global coupling strength in the system, and $K_j$($>0$) corresponds to the effective coupling strength to an oscillator $j$ from the population. All oscillators possess identical intrinsic frequency $\omega_j = \omega$. $\lambda_j$ is the bifurcation parameter controlling how fast the trajectory decays onto the attractor; in this model we consider $\lambda_j = \lambda > 0$ for all $j = 1, 2, 
3 ... N$, such that $\sqrt{\lambda}$ is considered as the "intrinsic amplitude" to which the oscillator converges in the absence of the coupling. $\beta$ is a phase delay term, and $d_0 e^{-i\alpha}$ is a constant that translates $z_j$ by a fixed amount in both amplitude and phase. (Alternatively, one could consider $\alpha \in [0,2\pi)$ and $d_0 \in \mathbb{R}>0$, but for simplicity in analysis we consider the above parameters). Note that Eq. (\ref{eq:SL1}) is a generalized form of diffusively coupled Stuart-Landau oscillators, such that when $\alpha = 0,\, \beta = 0$, and $d_{0} = 1$, the coupling function becomes $H = z_k - z_j$, which has been studied extensively.\cite{yamaguchi1984theory, bonilla1987self, matthews1990phase, matthews1991dynamics, monte2002dynamics, atay2003distributed, zakharova2014chimera} In polar coordinates, the model is written as
\begin{eqnarray} \label{eq:SL2}
\dot{\theta}_{j} &=& \omega + \frac{S K_j}{N}\sum_{k=1}^{N}[\frac{r_k}{r_j}\sin(\theta_k-\theta_j-\beta)+d_{0}\sin\alpha], \\ 
\dot{r}_{j} &=& (\lambda - |r_j|^{2})r_j \nonumber \\ 
&+& \frac{S K_j}{N}\sum_{k=1}^{N}[r_{k}\cos(\theta_k-\theta_j-\beta) -r_{j} d_{0}\cos\alpha], \label{eq:SL3}
\end{eqnarray}
where $\theta_{j}(t)$ is the phase and $r_{j}(t)$ is the amplitude of an oscillator $j$ at time $t$. 
This is a direct extension of the phase-reduced model investigated in Ref. \onlinecite{chaos2019}. The phase-reduced model is derived from Eq. (\ref{eq:SL2}) by setting the amplitudes of all oscillators to be constant and equal.  As explained further in Section \ref{sec:fullConnectivity}, the mean-field model serves as an approximation of a full network model with sufficiently large $N$, where $K_j$, the coupling strength, is directly proportional to the degree of node $j$.\cite{twko2008bistable, twko2008partial}

We investigate the effect of $\beta$, $\alpha$, $d_0$ and $\{K_j\}$ on the dynamics of the coupled oscillators. $\{K_j\}$ describes the distribution of coupling strengths, while different values of $\beta$, $\alpha$, and $d_0$ further determine the form of coupling function between nodes. Note that the main source of inhomogeneity in the system is through inhomogeneous values of $K_j$ for the oscillators. Instead of varying the intrinsic properties of oscillators through $\omega_j$, which has been an often adopted approach for many previous studies\cite{strogatz1991stability, bonilla1992nonlinear, crawford1994amplitude}, we focus on the dynamics created through inhomogeneous coupling strengths while assuming identical oscillators with $\omega_j = \omega$. This has a more practical implication for understanding real-world complex networks: in brain networks, for example, individual patches of cortical tissue are often modeled as identical units, differentiated only by the pattern and strength of connections they maintain with other nodes.\cite{sporns2014contributions} \newline

To study the effect of $d_0$, $\alpha$, $\beta$, and $\{K_j\}$ on the dynamics of the system, we perform a similar self-consistency analysis as in Ref. \onlinecite{chaos2019} with the added amplitude dimension. To simplify our notation, we set S = 1 in the following analysis without loss of generality.

Let $\Omega$ denote the frequency of the population oscillation described by the order parameter $\tilde{R} e^{i \Theta} \equiv \frac{1}{N} \sum_{j=1}^{N} r_j e^{i \theta_{j}}$ in a stationary state. 
Then
\begin{eqnarray}
\dot \phi_j &=& \Delta + K_j \left [d_0 \sin\alpha+ \frac{1}{r_j}\tilde{R} \sin(\Phi-\phi_j-\beta) \right],
\label{eq_sine_1}
\end{eqnarray}
where $\phi_j \equiv \theta_j - \Omega t$, $\Delta \equiv \omega-\Omega$, 
and $\Phi \equiv \Theta - \Omega t$.
When the system reaches a stationary state, $\tilde{R}$ and $\Phi$ do not depend on time.  
Additionally,
\begin{eqnarray}
\dot{r}_j &=& \lambda r_j - {r_j}^3 + K_j \left [\tilde{R}\cos(\Phi-\phi_j-\beta) - r_j d_0\cos\alpha \right]. \nonumber \\ 
\label{eq_cosine1}
\end{eqnarray}

The oscillators phase-locked with frequency $\Omega$ in the original frame of reference are those with $K_j \in \mathcal{D}_{l} \equiv \{K_j : K_j \tilde{R} > |\Delta + K_j d_0 \sin\alpha|{r_j}^* \:\}$ asymptotically approaching a stable fixed point $z^*_j$ = (${\phi_j}^*$, ${r_j}^*$) of Eqs. (\ref{eq_sine_1}) and (\ref{eq_cosine1}), satisfying  the following equations:
\begin{eqnarray}
(\Delta + K_j d_0\sin\alpha)r^*_j &=& K_j \tilde{R} \sin\left({\phi_j}^* - \Phi + \beta \right)
\label{eq_Plocked}	\\
(\lambda - {r^*_j}^2 - K_j d_0\cos\alpha)r^*_j &=& -K_j \tilde{R} \cos\left({\phi_j}^* - \Phi + \beta \right) \nonumber \\ 
\label{eq_Rlocked}
\end{eqnarray}
from $\dot \phi_j = 0$ and $\dot r_j = 0$ respectively. This coupled system yields an exact equation for ${r_j}^*$ as 
\begin{align}
\{ (\Delta + K_j d_0\sin\alpha)^2 + (\lambda - {r^*_j}^2 - K_j d_0\cos\alpha)^2 \}{r^*_j}^2 \nonumber \\= {(K_j\tilde{R})}^2 . 
\label{amp_eq}
\end{align}
Furthermore,
\begin{eqnarray}
\label{eq_Plocked_stabcond}
\cos\left({\phi_j}^* - \Phi + \beta \right) > 0	\\
\label{eq_Rlocked_stabcond}
\lambda - 3{r^*_j}^2 - K_j d_0 \cos \alpha < 0
\end{eqnarray}
due to the stability of the fixed point. Combining Eqs. (\ref{eq_Rlocked}) and (\ref{eq_Plocked_stabcond}) gives 
\begin{eqnarray}
\lambda - {r^*_j}^2 - K_j d_0 \cos \alpha < 0 
\label{eq_stabcond}
\end{eqnarray}
which, when applied to Eq. (\ref{amp_eq}), yields one real positive solution for $r^*_j$.

From Eqs. (\ref{eq_Plocked}), (\ref{eq_Rlocked})  and (\ref{eq_Plocked_stabcond}), we obtain the fixed points:
\begin{eqnarray} 
{\phi_j}^*  &=& \sin^{-1}\left[\frac{(\Delta + K_j d_0 \sin \alpha)r^*_j}{K_j \tilde{R}}\right]+\Phi - \beta 
\label{eq:locked_phases1} \\
&=& \cos^{-1}\left[\frac{-(\lambda -{r^*_j}^2 - K_j d_0 \cos \alpha)r^*_j}{K_j \tilde{R}}\right]+\Phi - \beta. \nonumber \\ 
\label{eq:locked_phases2}
\end{eqnarray}

Here we consider two curves $(K_j, {\phi_j}^*)$ and $(K_j, {r_j}^*)$, defined as the distributions of ${\phi_j}^*$ and ${r_j}^*$ as functions of $K_j$. The slopes of these curves can provide useful information for understanding the dynamics of the oscillators with regard to coupling strength. First, the slope of $(K_j, {\phi_j}^*)$ curve for locked oscillators describes the phase distribution over the oscillators with regard to the distribution of the coupling strength. This slope is given by
\begin{eqnarray} 
\frac{\partial {\phi_j}^*}{\partial K_j}
= -\frac{\Delta \cdot r^*_j}{{K_j}^2 \tilde{R} \cos\left({\phi_j}^* - \Phi + \beta \right) }.
\label{eq:sign_slope_phi}
\end{eqnarray}

With the condition Eq. (\ref{eq_Plocked_stabcond}), the sign of the slope further reduces to 
\begin{eqnarray} 
\rm{sign} \left (\frac{\partial {\phi_\textit{j}}^*}{\partial \textit{K}_\textit{j} }\right )
= -\rm{sign}\left(\Delta \cdot \textit{r}^*_\textit{j} \right)
= -\rm{sign}\left(\Delta\right)
\label{eq:sign_slope_phi2}
\end{eqnarray}
since $r^*_j > 0$ at the stable point. This means only the sign of $\Delta$ determines the sign of the slope and the phase monotonically increases or decreases as $K_j$ increases within the locking range of $K_j$ values. When the locked oscillators oscillate with frequency $\Omega$ greater than the intrinsic frequency $\omega$ ($\Delta < 0$), the slope is positive. On the other hand, the oscillators with $K_j \in \mathcal{D}_{d} \equiv \{K_j : K_j \tilde{R} < |\Delta + K_j d_0 \sin\alpha|r^*_j \:\}$ drift monotonically without locking. We refer to these oscillators as drifting population.

In order to reveal how the amplitudes of the locked oscillators change as a function of the coupling strength, slope of the $(K_j, r_j)$ curve can next be found. From Eq. (\ref{eq_Plocked}) we can derive
\begin{eqnarray} 
\frac{\partial {r_j}^*}{\partial K_j}
= \frac{\Delta \tilde{R}\sin\left (\phi^*_j-\Phi+\beta\right)}{\left(\Delta+K_j d_0 \sin\alpha\right)^2}.
\end{eqnarray}
which reduces to
\begin{eqnarray} 
\rm{sign} \left (\frac{\partial {\textit{r}_\textit{j}}^*}{\partial \textit{K}_\textit{j} }\right )
&=& \rm{sign}\left \{\Delta \sin\left( \phi^*_\textit{j}-\Phi+\beta \right) \right\} \\
&=&
\begin{cases}
\rm{sign}\left(\Delta\right) &\text{if $\phi^*_j-\Phi+\beta \in (0, \frac{\pi}{2})$} \\
-\rm{sign}\left(\Delta\right) &\text{if $\phi^*_j-\Phi+\beta \in (-\frac{\pi}{2}, 0)$} \\
0 &\text{if $\phi^*_j-\Phi+\beta = 0$} .
\end{cases} \nonumber
\label{eq:sign_slope_r}
\end{eqnarray}
This indicates that the slope is not necessarily monotonic. Since $\Phi = 0$ at the stationary state by definition, the values of $\phi^*$ and $\beta$ determine the sign of the slope with respect to $\Delta$. For example, if $\phi^*_j+ \beta$ falls in the range of $(-\frac{\pi}{2}, \frac{\pi}{2})$, then the inflection point will occur at where $\phi^*_j+ \beta = 0$.

We now calculate the order parameter contributions from the locked and drifting subpopulations in the rotating frame. The self-consistency condition requires that $\tilde{R} = \tilde{R}_{lock} + \tilde{R}_{drift}$. 
Noting that $e^{i(\phi - \Phi + \beta)} = \cos(\phi - \Phi + \beta) + i \sin(\phi - \Phi + \beta)$ and using Eq. (\ref{eq_Plocked}), the contribution from the locked oscillators can be calculated as follows for a given coupling strength distribution $g(K)$:
\begin{eqnarray} 
\tilde{R}_{l} &=& \int_{\mathcal{D}_{l}} g(K) r^*e^{i\phi} dK \nonumber \\
&=& e^{-i\beta} \int_{\mathcal{D}_{l}} g(K) r^* \frac{\sqrt{K^2{\tilde{R}}^2-(\Delta^{'})^2 {r^*}^2} + i \Delta^{'}}{K\tilde{R}} dK, \nonumber \\
\label{eq:R_lock}
\end{eqnarray}
where the abbreviation $\Delta^{'} \equiv \Delta + K d_0 \sin \alpha$   is used for simplification. Combining Eqs. (\ref{amp_eq}) and (\ref{eq_stabcond}) yields an exact solution for $r^*$.

In order to determine the order parameter contribution from the drifting subpopulation, we follow the similar perturbation method used in Ref. \onlinecite{matthews1991dynamics}. The difference here is that the current model replaces the term $\Delta$ with the added desynchronization factor $\Delta + K_{j}d_{0}\sin \alpha$ as shown in Eq. (\ref{eq_sine_1}). 
In incoherent state, $\tilde{R}=0$ and $r= \sqrt{\lambda-K d_0\cos\alpha} \equiv a$. To find the bifurcation, we give a small perturbation $\epsilon$ to the system such that $\tilde{R}=\epsilon\tilde{R_1}$ where $\tilde{R_1}$ is a nonzero constant, $r = a + O(\epsilon)$, and $\phi_j(t) = \phi_0 + (\Delta+K_j d_0\sin\alpha)t + O(\epsilon)$. We seek a self-consistent partially locked solution to the first order in $\epsilon$ by finding the bifurcation from the incoherent state. By perturbing around the $\epsilon=0$ solution $r_j(t) = a$, $\phi_j(t) = (\Delta+K_j d_0\sin\alpha)t + \phi_0$, we obtain the path of the limit cycle as
 $r_j=a + \epsilon [A \cos(\phi-\Phi+\beta) + B \sin(\phi-\Phi+\beta) ] + O(\epsilon^2)$, where $A = \frac{2a^2}{(\Delta^{'})^2+4a^4}$ and $B = \frac{\Delta^{'}}{(\Delta^{'})^2+4a^4}$. We now impose the stationary condition requiring that oscillators with coupling strength $K_j$ form a stationary distribution along their limit cycles. Then the contribution to the centroid from the drifting oscillator of degree $K_j$ is $R_d (K_j) = e^{-i\beta}\epsilon(A+iB)/2 + O(\epsilon^2)$. Therefore,
\begin{eqnarray} 
\tilde{R}_{d} &=& \int_{\mathcal{D}_{d}} g(K) R_d(K_j) dK + O(\epsilon^2) \nonumber \\
&=& \frac{1}{2} e^{-i\beta} \int_{\mathcal{D}_{d}} g(K) K\tilde{R} \frac{2a^2 + i\Delta^{'}}{\Delta^{'}+4a^4} dK + O(\epsilon^2).  \nonumber \\
\label{eq:R_drift}
\end{eqnarray}
Thus from $\tilde{R} = \tilde{R}_{l} + \tilde{R}_{d}$ we obtain two independent equations for the values of $\tilde{R}$ and $\Delta$, which we can solve numerically for given $\alpha$, $\beta$,  $d_0$, and $g(K)$. We note that other non-perturbative methods may be applicable.\cite{clusella2019between}

\begin{table*}
\centerfloat
\caption{Categorization of the synchronous states, for the case where $d_0\sin\alpha \geq 0$: The name of state is given as $Sn_x$, following the same naming scheme presented in Ref \onlinecite{chaos2019}: $n$ is the major category index, and $x$ is composed of $d$, $l^+$, $l^-$, and $l^0$ where $d$ stands for a drifting range of $K$, $l$ for a locking range of $K$, $l^{+}$, $l^{-}$, $l^{0}$, respectively, for positive slope, negative slope, and zero slope of the curve $(K_j, {\phi_j}^*)$ in the locking range of $K$. $\Delta \equiv \omega - \Omega$. 
$D_0 \equiv |d_0\sin\alpha|{r_j}^*$ in the last column with ${r_j}^*={r_{\rm min}}^*$ for $K_{\rm min}$ and ${r_j}^*={r_{\rm max}}^*$ for $K_{\rm max}$. For other details, see the text.}
\begin{tabular}
{>{\let\centering\relax}m{1cm}  >{\centering}m{3.5cm} >{\centering}m{1.5cm} >{\centering}m{1.5cm} >{\centering}m{1.5cm} >{\centering}m{3.5cm} >{\centering}m{4.7cm}}
\hline\hline 
States   & Oscillators with $K$ \\from $K_{\rm min}$  to  $K_{\rm max}$  & Slope of \\$(K_j, {\phi_j}^*)$ & Sign$(\Delta)$ & $(\tilde{R}, D_0)$ & Locking range of $K$ & Additional condition \tabularnewline
\hline
$S1_{l^0}$        & In-phase synchronous           & $0$ &  $0$           &  $\tilde{R}>D_0$ &  $[K_{\rm min}, K_{\rm max}]$  & $\text{max} \, \tilde{R}, \Delta = 0$  \tabularnewline
$S1_{l^+}$        & Fully locked           & $+$ &  $-$           &  $\tilde{R}\geq D_0$ &  $[K_{\rm min}, K_{\rm max}]$  & $\frac{|\Delta|r^*_j}{\tilde{R}+D_0} < K_{\rm min}$ \tabularnewline
$S1_{dl^+}$       & Drifting$-$locked    & $+$&  $-$           &  $\tilde{R}\geq D_0$ & $\frac{|\Delta|r^*_j}{\tilde{R}+D_0} < K_j$ & $K_{\rm min} \leq \frac{|\Delta|r^*_j}{\tilde{R}+D_0}$ \tabularnewline
\hline 
$S2_{l^-}$        & Fully locked           & $-$ &  $+$           &  $\tilde{R} > D_0$  &  $[K_{\rm min}, K_{\rm max}]$	&	$ \frac{\Delta r^*_j}{\tilde{R}-D_0} < K_{\rm min}$	\tabularnewline
$S2_{dl^-}$       & Drifting$-$locked    & $-$ &  $+$           &  $\tilde{R} >D_0$ & $\frac{\Delta r^*_j}{\tilde{R}-D_0} < K_j$ & $K_{\rm min} \leq \frac{\Delta r^*_j}{\tilde{R}-D_0}$	 \tabularnewline
$S2_{d}$       & Fully drifting   & $-$ &  $+$           &  $\tilde{R} >D_0$ & None & $K_{\rm max} \leq \frac{\Delta r^*_j}{\tilde{R}-D_0}$	 \tabularnewline
\hline 
$S3_{l^+}$        & Fully locked           & $+$	&  $-$           &  $\tilde{R} < D_0$ &  $[K_{\rm min}, K_{\rm max}]$	&	$\frac{|\Delta|r^*_j}{\tilde{R}+D_0} < K_{\rm min}$, $K_{\rm max}  < \frac{|\Delta|r^*_j}{D_0-\tilde{R}}$	\tabularnewline
$S3_{l^+d}$       & Locked$-$drifting    & $+$ &  $-$           &  $\tilde{R} < D_0$ & $K_j < \frac{|\Delta|r^*_j}{D_0-\tilde{R}}$ &	$\frac{|\Delta|r^*_j}{\tilde{R}+D_0} < K_{\rm min}$, $\frac{|\Delta|r^*_j}{D_0-\tilde{R}} \leq K_{\rm max}$\tabularnewline
$S3_{dl^+}$       & Drifting$-$locked    & $+$ &  $-$           &  $\tilde{R} < D_0$ & $\frac{|\Delta|r^*_j}{\tilde{R}+D_0} < K_j$ &	$K_{\rm min} \leq \frac{|\Delta|r^*_j}{\tilde{R}+D_0} $, $K_{\rm max}  < \frac{|\Delta|r^*_j}{D_0-\tilde{R}}$	\tabularnewline
$S3_{dl^+d}$      & Drifting$-$locked$-$drifting& $+$ &  $-$  &  $\tilde{R} < D_0$ & $\frac{|\Delta|r^*_j}{\tilde{R}+D_0} < K_j < \frac{|\Delta|r^*_j}{D_0-\tilde{R}}$	&	$K_{\rm min} \leq \frac{|\Delta|r^*_j}{\tilde{R}+D_0}$, $\frac{|\Delta|r^*_j}{D_0-\tilde{R}} \leq K_{\rm max} $ \tabularnewline
$S3_{d}$          & Fully drifting      & None &  	$-$             & $\tilde{R} < D_0$ & None  & $\frac{|\Delta|r^*_j}{D_0-\tilde{R}} \leq K_{\rm min}$ or $K_{\rm max} \leq \frac{|\Delta|r^*_j}{\tilde{R}+D_0}$ \tabularnewline
\hline 
$S4_{d}$          & Fully drifting      & None &  $+,0$             & $\tilde{R} \leq D_0$ & None  & ... \tabularnewline
\hline\hline
\end{tabular}
\label{table:table_1}
\end{table*}

\begin{table*}
\centerfloat
\caption{Categorization of the synchronous states, for the case where $d_0\sin\alpha < 0$. $D_0 \equiv |d_0\sin\alpha|{r_j}^*$ in the last column with ${r_j}^*={r_{\rm min}}^*$ for $K_{\rm min}$ and ${r_j}^*={r_{\rm max}}^*$ for $K_{\rm max}$. Other details are as in Table I. Note the difference in signs of $\Delta$ and the boundary values in the additional conditions.}
\begin{tabular}
{>{\let\centering\relax}m{1cm}  >{\centering}m{3.5cm} >{\centering}m{1.5cm} >{\centering}m{1.5cm} >{\centering}m{1.5cm} >{\centering}m{3.5cm} >{\centering}m{4.7cm}}
\hline\hline 
States   & Oscillators with $K$ \\from $K_{\rm min}$  to  $K_{\rm max}$  & Slope of \\$(K_j, {\phi_j}^*)$ &  Sign$(\Delta)$ & $ (\tilde{R}, D_0)$ & Locking range of $K$ & Additional condition \tabularnewline
\hline
$S1_{l^0}$        & In-phase synchronous           & $0$ &  $0$           &  $\tilde{R}>D_0$ &  $[K_{\rm min}, K_{\rm max}]$  & $\text{max} \, \tilde{R}, \Delta = 0$  \tabularnewline
$S1_{l^+}$        & Fully locked           & $+$ &  $-$           &  $\tilde{R} > D_0$ &  $[K_{\rm min}, K_{\rm max}]$  & $\frac{|\Delta|r^*_j}{\tilde{R}+D_0} < K_{\rm min}$ \tabularnewline
$S1_{dl^+}$       & Drifting$-$locked    & $+$&  $-$           &  $\tilde{R} > D_0$ & $\frac{|\Delta|r^*_j}{\tilde{R}+D_0} < K_j$ & $K_{\rm min} \leq \frac{|\Delta|r^*_j}{\tilde{R}+D_0}$ \tabularnewline
\hline 
$S2_{l^-}$        & Fully locked           & $-$ &  $+$           &  $\tilde{R} \geq D_0$  &  $[K_{\rm min}, K_{\rm max}]$	&	$ \frac{\Delta r^*_j}{\tilde{R}+D_0} < K_{\rm min}$	\tabularnewline
$S2_{dl^-}$       & Drifting$-$locked    & $-$ & $+$           &  $\tilde{R} \geq D_0$ & $\frac{\Delta r^*_j}{\tilde{R}+D_0} < K_j$ & $K_{\rm min} \leq \frac{\Delta r^*_j}{\tilde{R}+D_0}$	 \tabularnewline
$S2_{d}$       & Fully drifting   & $-$ &  $+$           &  $\tilde{R} \geq D_0$ & None & $K_{\rm max} \leq \frac{\Delta r^*_j}{\tilde{R}+D_0}$	 \tabularnewline
\hline 
$S3_{l^-}$        & Fully locked           & $+$	&  $+$           &  $\tilde{R} < D_0$ &  $[K_{\rm min}, K_{\rm max}]$	&	$\frac{\Delta r^*_j}{\tilde{R}+D_0} < K_{\rm min}$, $K_{\rm max}  < \frac{\Delta r^*_j}{|\tilde{R}-D_0|}$	\tabularnewline
$S3_{l^-d}$       & Locked$-$drifting    & $+$ &  $+$           &  $\tilde{R} < D_0$ & $K_j < \frac{|\Delta|r^*_j}{|\tilde{R}-D_0|}$ &	$\frac{\Delta r^*_j}{\tilde{R}+D_0} < K_{\rm min}$, $\frac{\Delta r^*_j}{|\tilde{R}-D_0|} \leq K_{\rm max}$\tabularnewline
$S3_{dl^-}$       & Drifting$-$locked    & $+$ &$+$           &  $\tilde{R} < D_0$ & $\frac{\Delta r^*_j}{\tilde{R}+D_0} < K_j$ &	$K_{\rm min} \leq \frac{\Delta r^*_j}{\tilde{R}+D_0} $, $K_{\rm max}  < \frac{\Delta r^*_j}{|\tilde{R}-D_0|}$	\tabularnewline
$S3_{dl^-d}$      & Drifting$-$locked$-$drifting& $+$ & $+$  &  $\tilde{R} < D_0$ & $\frac{\Delta r^*_j}{\tilde{R}+D_0} < K_j < \frac{\Delta r^*_j}{|\tilde{R}-D_0|}$	&	$K_{\rm min} \leq \frac{\Delta r^*_j}{\tilde{R}+D_0}$, $\frac{\Delta r^*_j}{|\tilde{R}-D_0|} \leq K_{\rm max} $ \tabularnewline
$S3_{d}$          & Fully drifting      & None &  	$+$             & $\tilde{R} < D_0$ & None  & $\frac{\Delta r^*_j}{|\tilde{R}-D_0|} \leq K_{\rm min}$ or $K_{\rm max} \leq \frac{\Delta r^*_j}{\tilde{R}+D_0}$ \tabularnewline
\hline 
$S4_{d}$          & Fully drifting      & None & 	$0, -$             & $\tilde{R} \leq D_0$ & None  & ... \tabularnewline
\hline\hline
\end{tabular}
\label{table:table_2}
\end{table*}

Now, let us determine the range of $K_j$ for which the oscillators are phase-locked and classify the locked states. We can find the "locking" range of $K_{j}$ by jointly considering the sign of $\sin\alpha$, the sign of $\Delta + K_j d_0 \sin\alpha$, the sign of $\Delta$, and the locking condition $K_j \tilde{R} > |\Delta + K_j d_0 \sin\alpha|r^*_j $. The range of $K_j$ and the slope of $(K_j, {\phi_j}^*)$ for the locked oscillators are as follows: \\
\\ (i) If $d_0\sin\alpha \geq 0$:
\begin{align}
\begin{cases}
(a) ~S1: \frac{|\Delta|r^*_j}{\tilde{R}+d_0\sin\alpha\, r^*_j} < K_j, ~l^+ ~~~
\text{if $\tilde{R} \geq d_0\sin\alpha\,r^*_j$, $\Delta \leq 0$} \\
(b) ~S2:\frac{\Delta r^*_j}{R-d_0\sin\alpha\, r^*_j} < K_j, ~l^-  ~~~ 
\text{if $\tilde{R} > d_0\sin\alpha\, r^*_j$, $\Delta > 0$} \\ 
(c) ~S3: \frac{|\Delta|r^*_j}{\tilde{R}+d_0\sin\alpha\, r^*_j} < K_j < \frac{|\Delta|r^*_j}{d_0\sin\alpha\, r^*_j-\tilde{R}}, ~l^+
\\ \kern 11.5pc \text{if $\tilde{R} < d_0\sin\alpha\, r^*_j$, $\Delta < 0$} \\
(d) ~S4: \text{no locking range}\kern 2.5pc \text{if $\tilde{R} < d_0\sin\alpha\, r^*_j$, $\Delta \geq 0$}
\label{pstate_conditions1}
\end{cases} 
\end{align}
\\ (ii) If $d_0\sin\alpha < 0$:
\begin{align}
\begin{cases}
(a) ~S1: \frac{|\Delta|r^*_j}{\tilde{R}+|d_0\sin\alpha| r^*_j} < K_j, ~l^+ ~~~
\text{if $\tilde{R} > |d_0\sin\alpha|r^*_j$, $\Delta \leq 0$} \\
(b) ~S2:\frac{\Delta r^*_j}{R-|d_0\sin\alpha| r^*_j} < K_j, ~l^-  ~~~ 
\text{if $\tilde{R} \geq |d_0\sin\alpha| r^*_j$, $\Delta > 0$} \\ 
(c) ~S3: \frac{\Delta r^*_j}{\tilde{R}+|d_0\sin\alpha|r^*_j} < K_j < \frac{\Delta r^*_j}{|d_0\sin\alpha| r^*_j-\tilde{R}}, ~l^-
\\ \kern 11.5pc \text{if $\tilde{R} < |d_0\sin\alpha| r^*_j$, $\Delta > 0$} \\
(d) ~S4: \text{no locking range}\kern 2.5pc \text{if $\tilde{R} \leq |d_0\sin\alpha| r^*_j$, $\Delta \leq 0$} 
\label{pstate_conditions2}
\end{cases} 	
\end{align}
where the two equalities of (ia) and (iid) do not hold at the same time. $l^-$, $l^0$, and $l^+$ represent the negative, the zero, and the positive slope of $(K_j, {\phi_j}^*)$ respectively. Overall, this analysis reveals the same set of states $S1$-$S4$ originally reported for the phase-reduced model in Ref \onlinecite{chaos2019}. However, one crucial difference from the previous work is that the set of states are now generalized to include amplitude dynamics that vary independently of the phase dimension. Furthermore, we have identified more general cases, differentiated by the sign of the constant term $d_0\sin\alpha$. In particular, the dynamics for $d_0\sin\alpha < 0$ demonstrate the possibility of new sub-states $S3_{l-}$, $S3_{l-d}$, $S3_{dl-}$, and $S3_{dl-d}$. Qualitatively, the S3 states are distinct from other states, because the high-degree nodes have the possibility to drift in phase. The conditions and the characteristics of these states are summarized in Table \ref{table:table_1} and \ref{table:table_2}. Table \ref{table:table_1} describes the states under condition (i), and Table \ref{table:table_2} describes those under condition (ii) above. Note that the sub-states are defined only in terms of dynamics in the phase dimension (e.g., locked/drifting). The amplitude dynamics can potentially differ within the same sub-state, as will be shown in the phase diagram in the next section.

\begin{figure}
\centering
\epsfig{figure=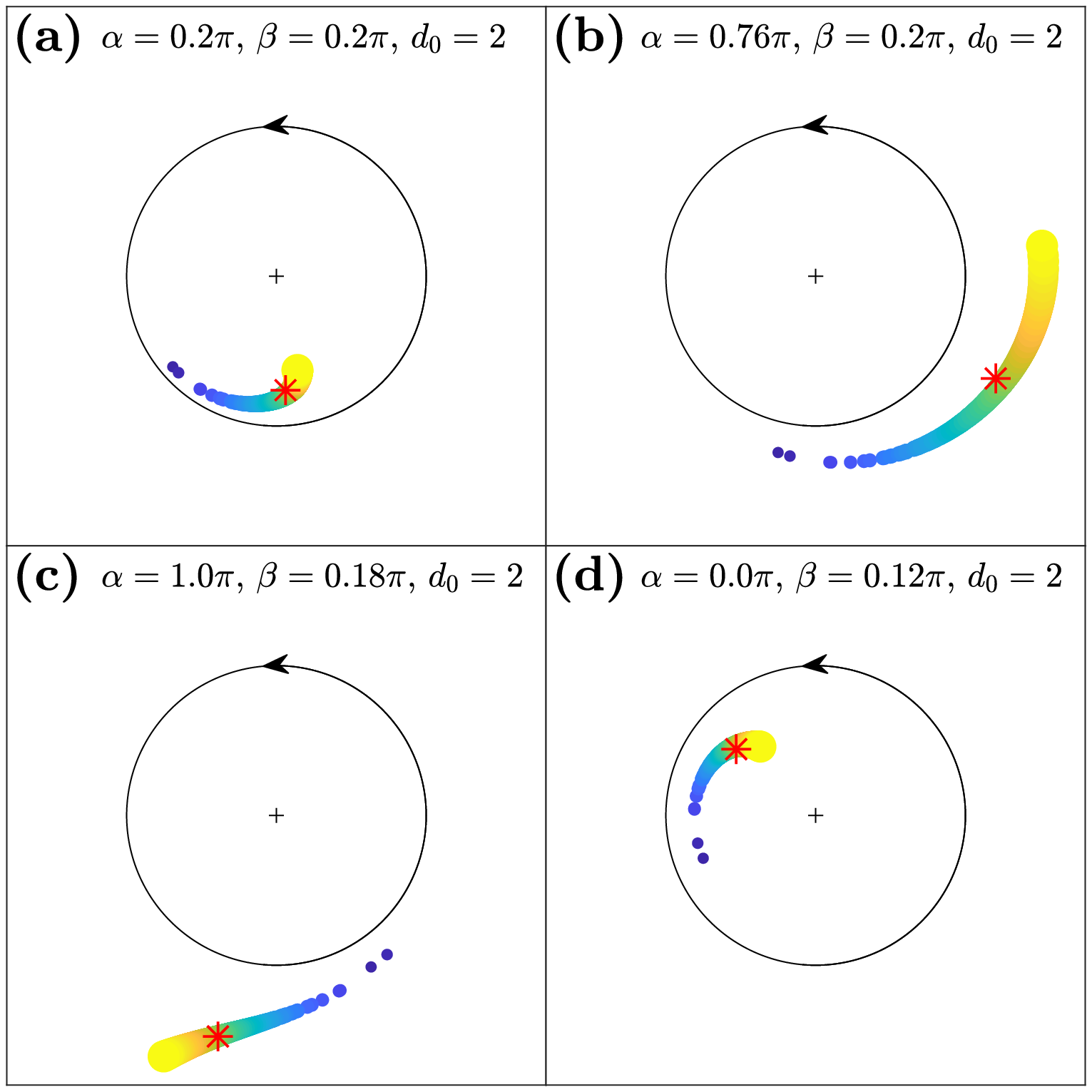, width=9cm}
\caption{ Synchronization of limit-cycle oscillators with inhomogeneous coupling strengths distribution (Gaussian). The amplitude and phase of each oscillator are represented geometrically on a complex plane. Color and size indicate the respective coupling strength---oscillators with larger coupling strengths are represented with bigger and brighter dots, and those with smaller coupling strengths are represented by smaller and darker dots. All oscillators are pulled towards the collective mean field (red asterisk) proportionately to each coupling strengths. The four states above show the representative synchronization dynamics for fully locked states. (a) Higher-degree nodes (with larger coupling strengths) phase-leading with lower amplitudes. (b) Higher-degree nodes phase-leading with higher amplitudes. (c) Higher-degree nodes phase-lagging with higher amplitudes. (d) Higher-degree nodes phase-lagging with lower amplitudes. 
Simulation parameters: $\lambda = 1$, $\omega = \pi$, $S = 12$, $N = 1000$. (Note that for this particular simulation, the global coupling strength $S$ was set to $12$ for the purpose of visualization; when $S$ is larger, the amplitude difference is exaggerated without a qualitative change in the dynamics; refer to Fig. S1 for more details.)
}
\label{fig1_ComplexPlane}
\end{figure}

\begin{figure*}
\centering
\epsfig{figure=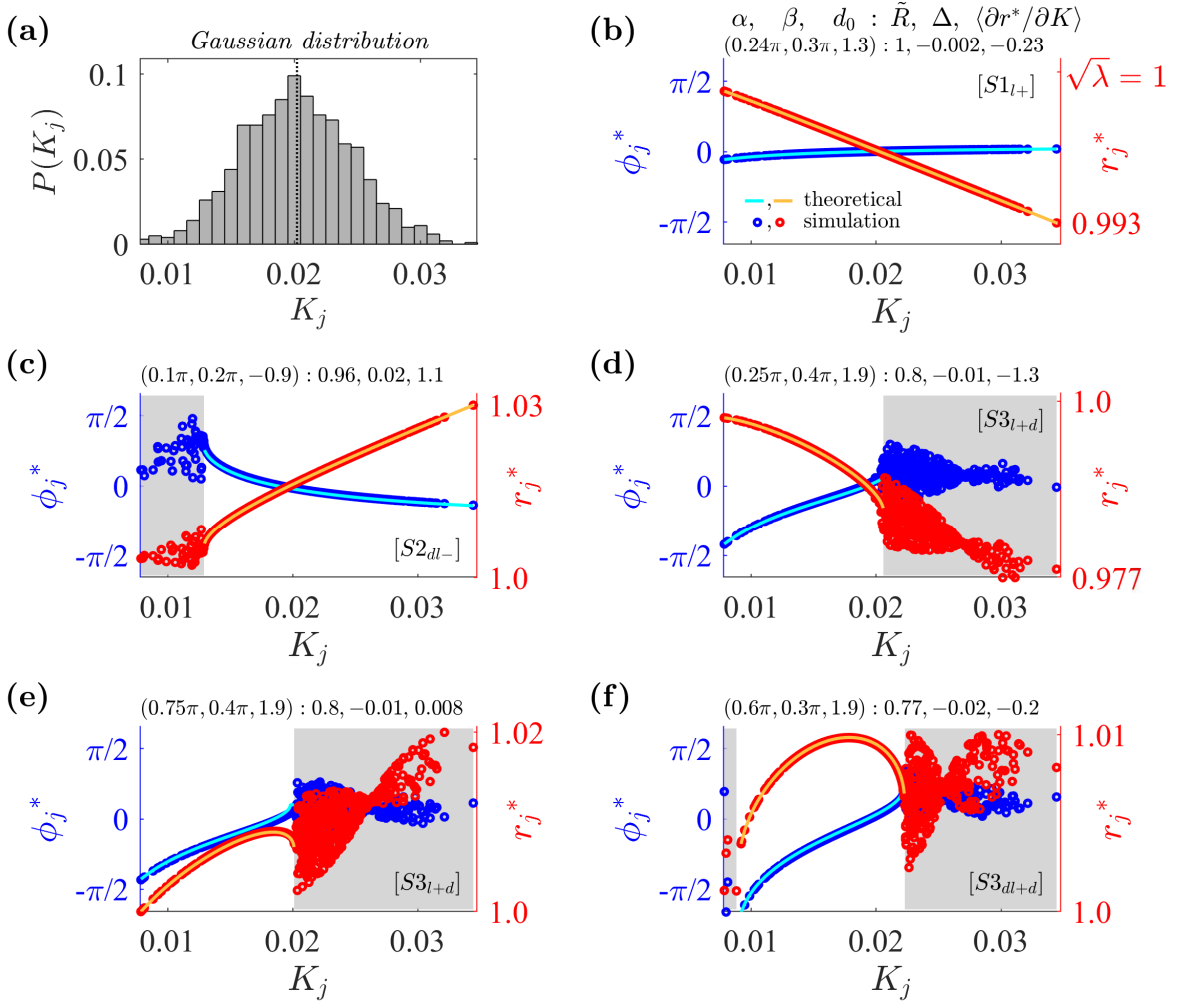, width=16cm}
\caption{Examples of various synchronous states in the system with Gaussian coupling strength distribution.
(a) Distribution of coupling strengths randomly selected from a Gaussian distribution with a mean of $20\times 10^{-3}$ and a standard deviation of $4.5\times 10^{-3}$. For the obtained coupling strength set ${K_j}$, $K_{\rm mean}=20.2\times 10^{-3}$ (vertical dotted line), $\sigma_K=4.36\times 10^{-3}$, $K_{\rm min} = 7.85\times 10^{-3}$, and $K_{\rm max} = 34.4\times 10^{-3}$.  
(b) State $S1_{l^+}$ with negative amplitude slope, where $l$ stands for the locking range of $K$: $\tilde{R} \geq d_0\sin\alpha\cdot{\rm max}\{{r_j}^*\}$ and  $\Delta <0$. Oscillators are fully locked with $\frac{|\Delta|{r_{\rm min}}^*}{\tilde{R}+|d_0\sin \alpha|{r_{\rm min}}^*} < K_{\rm min}$.
Left axis represents $\phi_j$ (in blue), and the right axis $r_j$ (in red).
(c) State $S2_{dl^-}$ with positive amplitude slope, where $d$ stands for the drifting range of $K$: $\tilde{R} \geq |d_0\sin\alpha|{\rm max}\{{r_j}^*\}$ and $\Delta >0$.  Oscillators with smaller coupling strengths drift with $K_{\rm min} \leq \frac{\Delta\cdot{r_{\rm min}}^*}{\tilde{R}+|d_0\sin \alpha|{r_{\rm min}}^*} < K_{max}$.
(d) State $S3_{l^+d}$ with negative amplitude slope: $\tilde{R} < d_0\sin\alpha\cdot{\rm min}\{{r_j}^*\}$ and $\Delta < 0$. Oscillators with larger coupling strengths drift with $\frac{|\Delta|{r_{\rm min}}^*}{\tilde{R}+ |d_0\sin \alpha|{r_{\rm min}}^*}< K_{\rm min}$ and $\frac{|\Delta|{r_{\rm max}}^*}{|d_0\sin \alpha|{r_{\rm max}}^*-\tilde{R}} < K_{\rm max}$. 
(e) State $S3_{l^+d}$ with positive-negative amplitude slope. Same additional conditions as in (d).
(f) State $S3_{dl^+d}$ with positive-negative amplitude slope: $\tilde{R} <  d_0\sin\alpha\cdot{\rm min}\{{r_j}^*\}$ and $\Delta < 0$. $K_{\rm min} \leq \frac{|\Delta|{r_{\rm min}}^*}{\tilde{R}+|d_0\sin \alpha|{r_{\rm min}}^*} < K_l < \frac{|\Delta|{r_{\rm max}}^*}{|d_0\sin \alpha|{r_{\rm max}}^*-\tilde{R}} \leq K_{\rm max}$, where $K_l$ represents the coupling strength for the locked oscillators.
In the figures (b)-(f), solid lines are self-consistent theoretical curves for locked phases and amplitudes from Eqs. (\ref{amp_eq}) and (\ref{eq:locked_phases1}), and unshaded range is for $K_j$ values for locked subpopulations obtained theoretically from Eqs. (\ref{pstate_conditions1}) and (\ref{pstate_conditions2}). 
Simulation parameters: $\lambda = 1$, $\omega = \pi$, $S = 1$, $N = 1000$.
}
\label{fig2_Gaussian_states}
\end{figure*}
\begin{figure*}
\centerfloat
\epsfig{figure=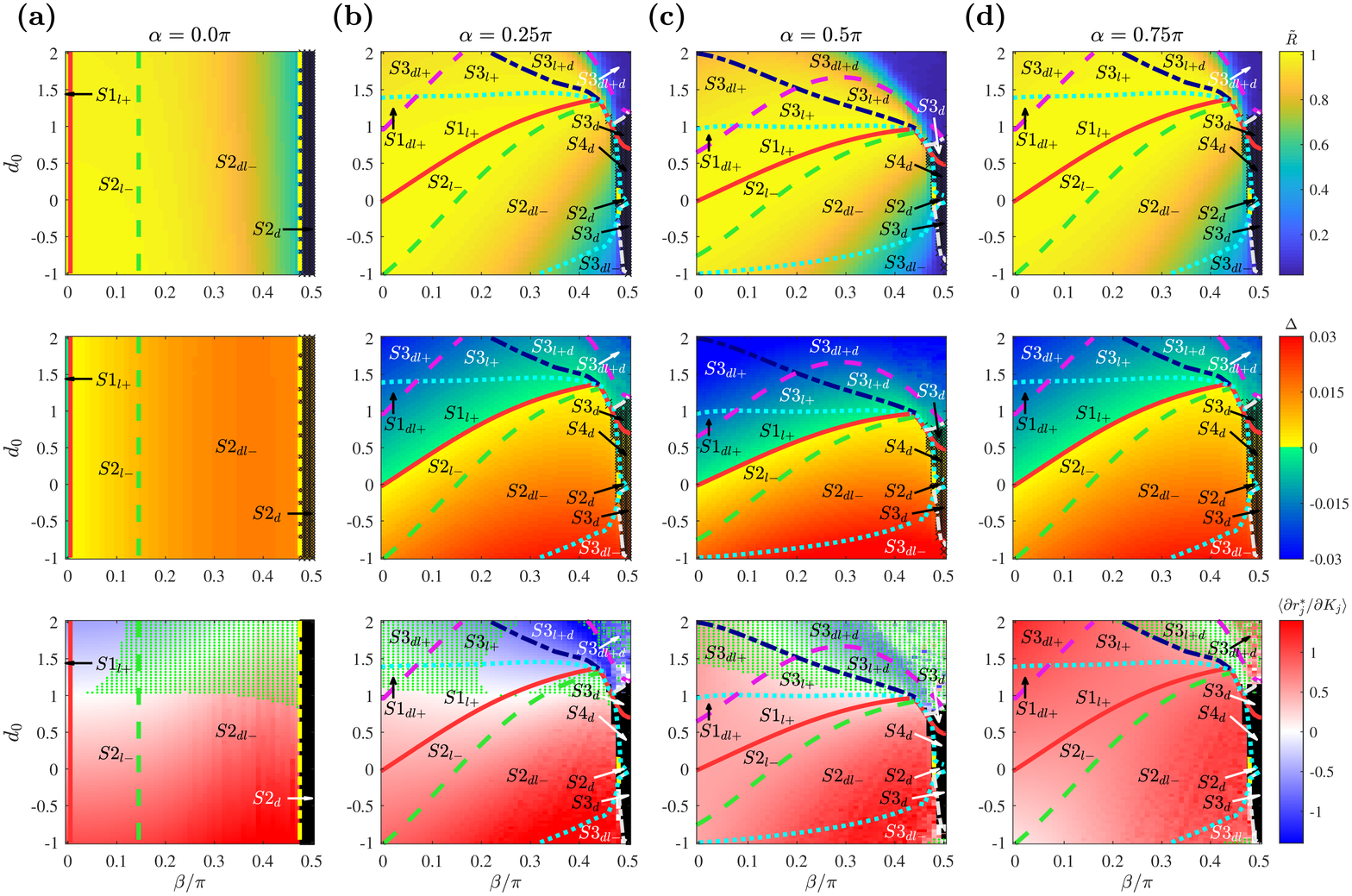, width=18.2cm}
\caption{
Phase diagrams with Gaussian coupling strength distribution as a function of $\alpha$, $\beta$, and $d_0$ determining the shape of coupling function.
Representative fixed values of $\alpha$ are chosen at (a) $\alpha = 0$ where $\sin\alpha = 0,\, \cos\alpha = 1$, (b) $\alpha = 0.25\pi$ where $\sin\alpha > 0,\, \cos\alpha > 0$, (c) $\alpha = 0.5\pi$ where $\sin\alpha = 1,\, \cos\alpha = 0$, and (d) $\alpha = 0.75\pi$ where $\sin\alpha < 0,\, \cos\alpha > 0$.
Phase diagrams with order parameter $\tilde{R}$ (first row), phase spread $\Delta = \omega - \Omega$ (second row), and slope of the amplitude curve $(K_j, {r_j}^*)$ calculated as the average value of $\frac{\partial r^*_j}{\partial K_j}$ among the locked oscillators (third row). Regions with fully drifting population are marked in black. For the amplitude slope, the parameter space in which inflection point exists in the $(K_j, {r_j}^*)$ curve (from positive to negative slope) is additionally marked with green dots. 
The boundaries are determined numerically from the model equation Eq. (\ref{eq:SL1}) for a Gaussian distribution with  $(K_{\rm mean}, \sigma_K, K_{\rm min}, K_{\rm max}) = (20.2,4.36,7.85,34.4) \times 10^{-3}$.
The bounding curves are obtained for $\Delta = 0$ (solid red), $\tilde{R} = |d_0\sin\alpha|\langle{r_j}^*\rangle$ (dotted cyan), $\frac{\Delta {r_{\rm min}}^*}{\tilde{R}-|d_0\sin\alpha|{r_{\rm min}}^*} = K_{min}$ (dashed green), $\frac{|\Delta|{r_{\rm min}}^*}{\tilde{R}+|d_0\sin\alpha|{r_{\rm min}}^*} = K_{min}$ (long dashed magenta), $\frac{|\Delta|{r_{\rm max}}^*}{|d_0\sin\alpha|{r_{\rm max}}^*-\tilde{R}} = K_{max}$ (dashed-dotted navy), $\frac{\Delta {r_{\rm max}}^*}{\tilde{R}-|d_0\sin\alpha|{r_{\rm max}}^*} = K_{max}$ (dashed-dotted yellow), $\frac{|\Delta|\sqrt{\lambda}}{|d_0\sin\alpha|\sqrt{\lambda}-\tilde{R}} = K_{min}$ (dashed-dotted grey, $d_0 \geq 0$), and $\frac{\Delta\sqrt{\lambda}}{\tilde{R}+|d_0\sin\alpha|\sqrt{\lambda}} = K_{max}$ (dashed-dotted grey, $d_0 < 0$). Simulation parameters: $\lambda = 1$, $\omega = \pi$, $S = 1$, $N = 1000$.
}
\label{fig3_Gaussian_diagram}
\end{figure*}
\begin{figure*}
\centering
\epsfig{figure=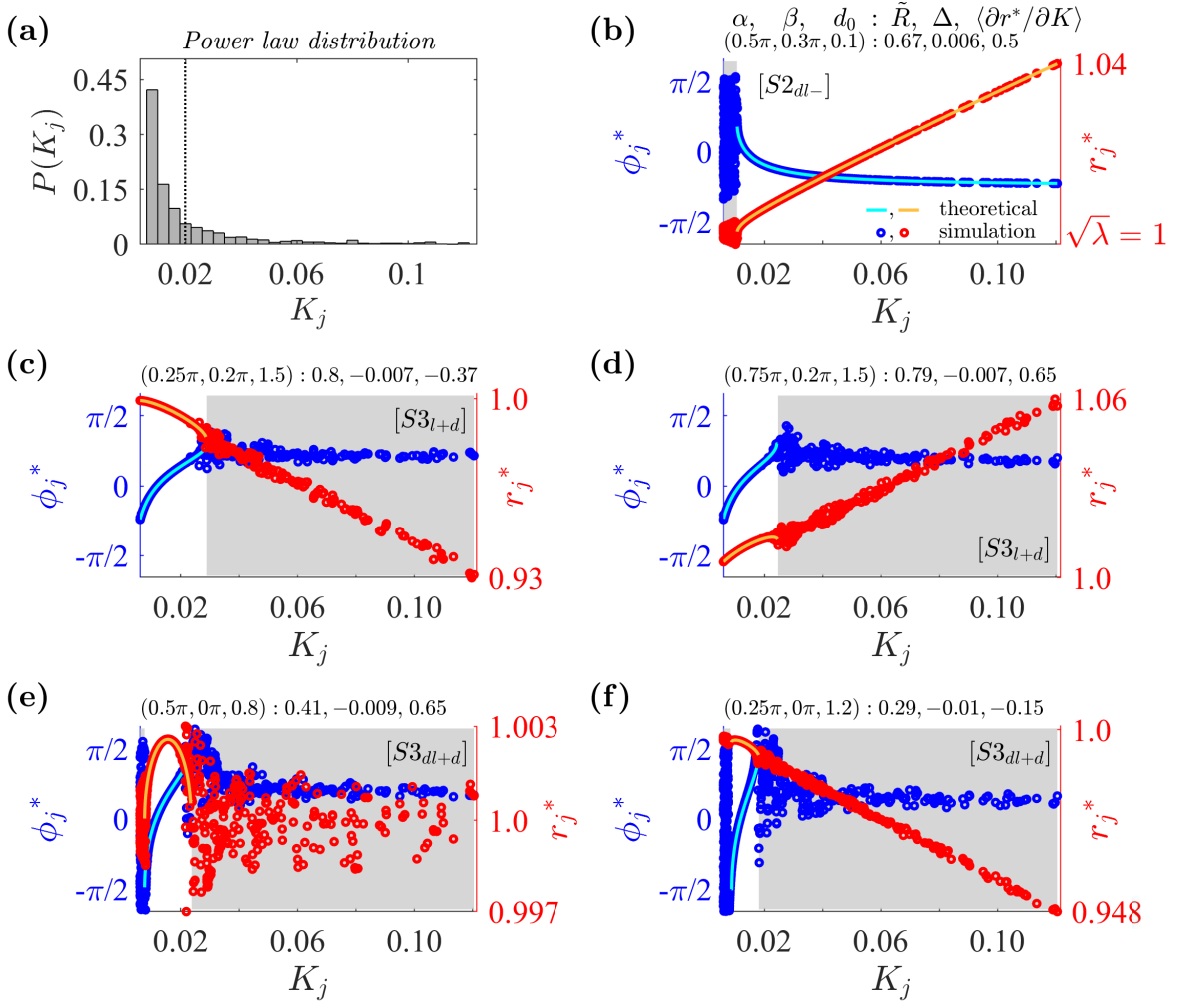, width=16cm}
\caption{Examples of synchronous states in the system with Power-law coupling strength distribution 
(a) Distribution for coupling strengths randomly selected from a power-law distribution with a mean of  $20\times 10^{-3}$ and an exponent of $2$. For the obtained coupling strength set ${K_j}$, $K_{\rm mean}=19.9 \times 10^{-3}$ (vertical dotted line) , $\sigma_K=20.5 \times 10^{-3}$, $K_{\rm min} =6.02\times 10^{-3}$, and $K_{\rm max} = 120.5 \times 10^{-3}$.  
(b) State $S2_{dl^-}$ with positive amplitude slope: $\tilde{R} > d_0\sin\alpha\cdot{\rm max}\{{r_j}^*\}$ and $\Delta > 0$. Satisfies additional condition $K_{\rm min} \leq \frac{\Delta\cdot{r_{\rm min}}^*}{\tilde{R}-|d_0\sin \alpha|{r_{\rm min}}^*} < K_{max}$. 
(c) State $S3_{l^+d}$ with negative amplitude slope: $\tilde{R} < D_0$ and $\Delta < 0$. Satisfies additional conditions $\frac{|\Delta|{r_{\rm min}}^*}{\tilde{R}+ |d_0\sin \alpha|{r_{\rm min}}^*}< K_{\rm min}$ and $\frac{|\Delta|{r_{\rm max}}^*}{|d_0\sin \alpha|{r_{\rm max}}^*-\tilde{R}} < K_{\rm max}$. 
(d) State $S3_{l^+d}$ with positive amplitude slope. Additional conditions are as in (c). Oscillators with larger coupling strengths drift with higher amplitudes.
(e) State $S3_{dl^+d}$ with positive-negative amplitude slope, $\tilde{R} <  d_0\sin\alpha\cdot{\rm min}\{{r_j}^*\}$ and $\Delta < 0$. Satisfies additional conditions
$K_{\rm min} \leq \frac{|\Delta|{r_{\rm min}}^*}{\tilde{R}+|d_0\sin \alpha|{r_{\rm min}}^*} < K_l < \frac{|\Delta|{r_{\rm max}}^*}{|d_0\sin \alpha|{r_{\rm max}}^*-\tilde{R}} \leq K_{\rm max}$. 
(f) State $S3_{dl^+d}$ with negative amplitude slope. Additional conditions are as in (e). 
Other details are as in Fig. \ref{fig2_Gaussian_states}. 
Simulation parameters: $\lambda = 1$, $\omega = \pi$, $S = 1$, $N = 1000$.
}
\label{fig4_SF_states}
\end{figure*}
\begin{figure*}
\centerfloat
\epsfig{figure=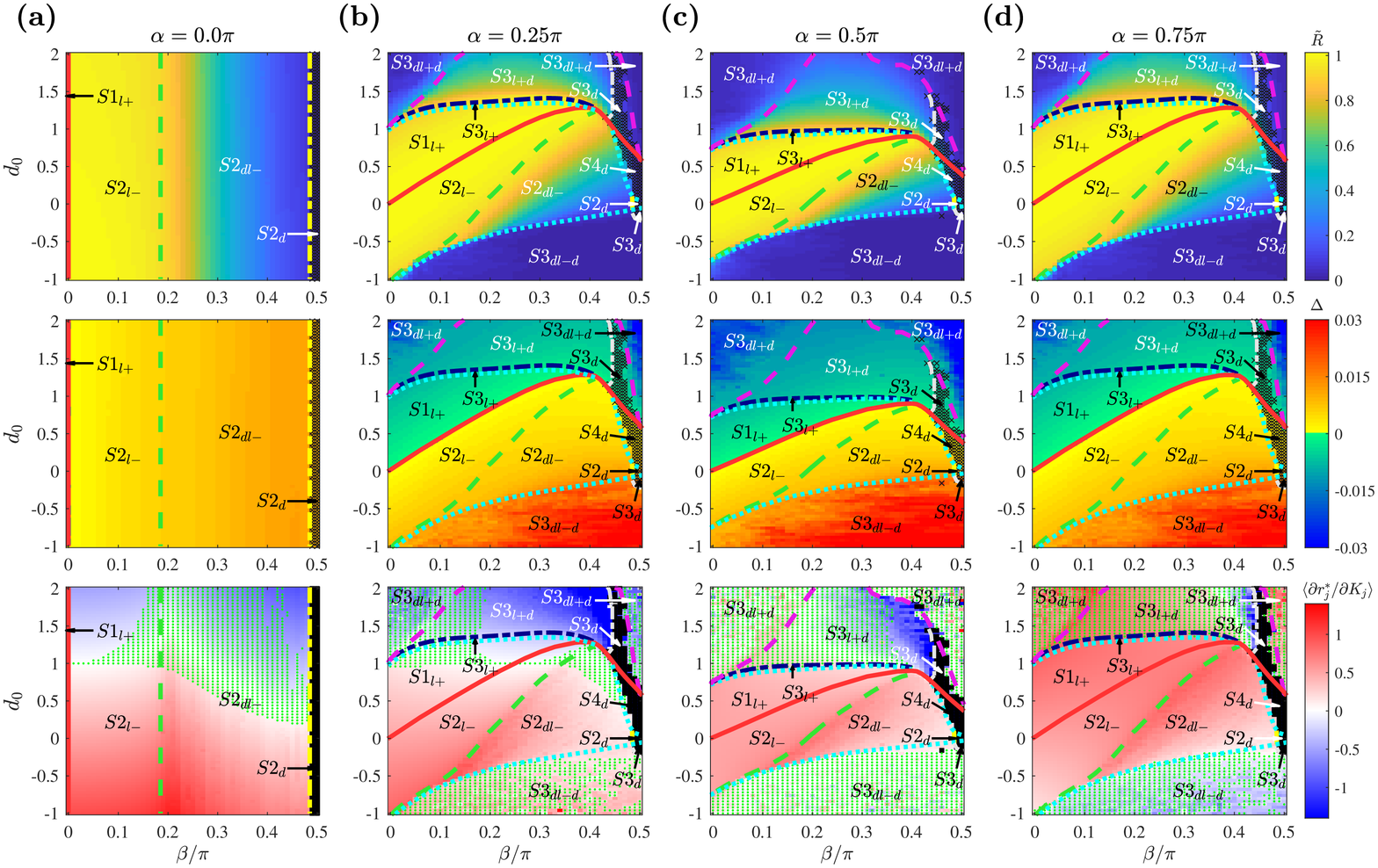, width=18.2cm}
\caption{
Phase diagrams with power-law coupling strength distribution as a function of $\alpha$, $\beta$, and $d_0$ determining the shape of coupling function. 
Representative fixed values of $\alpha$ are chosen at (a) $\alpha = 0$ where $\sin\alpha = 0,\, \cos\alpha = 1$, (b) $\alpha = 0.25\pi$ where $\sin\alpha > 0,\, \cos\alpha > 0$, (c) $\alpha = 0.5\pi$ where $\sin\alpha = 1,\, \cos\alpha = 0$, and (d) $\alpha = 0.75\pi$ where $\sin\alpha < 0,\, \cos\alpha > 0$.
Phase diagrams with order parameter $\tilde{R}$ (first row), phase spread $\Delta = \omega - \Omega$ (second row), and slope of the amplitude curve $(K_j, {r_j}^*)$ calculated as the average value of $\frac{\partial r^*_j}{\partial K_j}$ among the locked oscillators (third row). 
The boundaries and other details are as in Fig. \ref{fig3_Gaussian_diagram}.
The coupling strengths are randomly selected from a power-law distribution $P(x) \sim x^{-\gamma_0}$ with $\gamma_0 = 2$.  $\bar \gamma = 2.03$ is the average exponent of distributions for the obtained coupling strength sets. 
For the obtained coupling strength set $K_j$, $(K_{\rm mean}, \sigma_K, K_{\rm min}, K_{\rm max}) = (19.8, 20.2, 6.01, 122.7) \times 10^{-3}$. 
Simulation parameters: $\lambda = 1$, $\omega = \pi$, $S = 1$, $N = 1000$.
}
\label{fig5_SF_diagram}
\end{figure*}
\begin{figure*}
\centering
\epsfig{figure=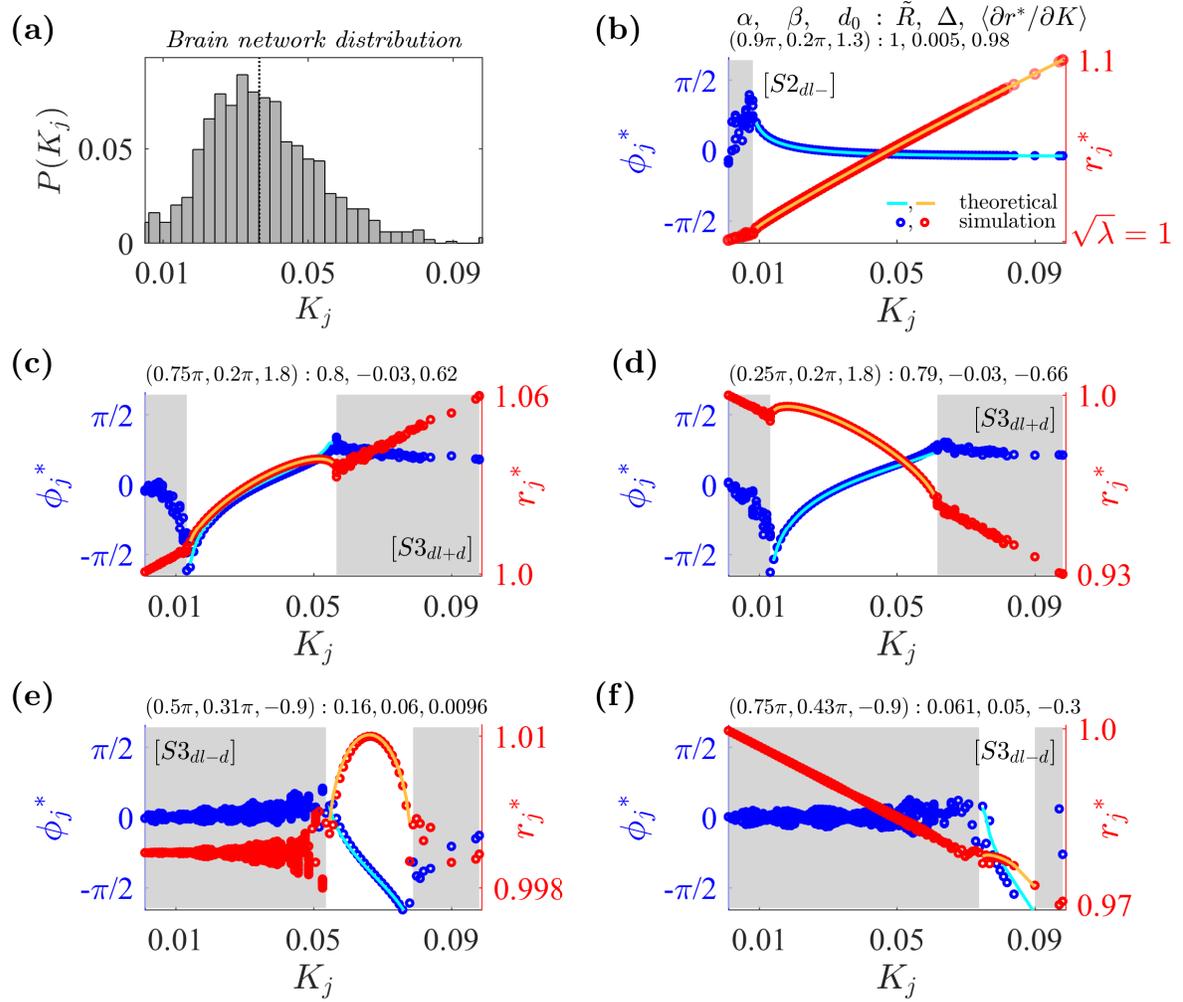, width=16cm}
\caption{Examples of various synchronous states in the system with human brain network coupling strength distribution. 
(a) Distribution for coupling strengths obtained from a brain network distribution. For the obtained coupling strength set ${K_j}$, $K_{\rm mean}=36.5 \times 10^{-3}$ (vertical dotted line), $\sigma_K=15.8 \times 10^{-3}$, $K_{\rm min} =1.01\times 10^{-3}$, and $K_{\rm max} = 98.1 \times 10^{-3}$.
(b) State $S2_{dl-}$ with positive amplitude slope: $\tilde{R} > d_0\sin\alpha\cdot{\rm max}\{{r_j}^*\}$ and $\Delta > 0$. Satisfies additional condition $K_{\rm min} \leq \frac{\Delta\cdot{r_{\rm min}}^*}{\tilde{R}-|d_0\sin \alpha|{r_{\rm min}}^*} < K_{max}$. 
(c) State $S3_{dl+d}$ with positive amplitude slope: $\tilde{R} <  d_0\sin\alpha\cdot{\rm min}\{{r_j}^*\}$ and $\Delta < 0$. Satisfies additional conditions
$K_{\rm min} \leq \frac{|\Delta|{r_{\rm min}}^*}{\tilde{R}+|d_0\sin \alpha|{r_{\rm min}}^*} < K_l < \frac{|\Delta|{r_{\rm max}}^*}{|d_0\sin \alpha|{r_{\rm max}}^*-\tilde{R}} \leq K_{\rm max}$. 
(d) State $S3_{dl+d}$ with negative amplitude slope. Additional conditions are as in (c).
(e) State $S3_{dl-d}$ with positive-negative amplitude slope: $\tilde{R} <  d_0\sin\alpha\cdot{\rm min}\{{r_j}^*\}$ and $\Delta > 0$. Satisfies additional conditions
$K_{\rm min} \leq \frac{|\Delta|{r_{\rm min}}^*}{\tilde{R}+|d_0\sin \alpha|{r_{\rm min}}^*} < K_l < \frac{|\Delta|{r_{\rm max}}^*}{|\tilde{R}-|d_0\sin \alpha|{r_{\rm max}}^*|} \leq K_{\rm max}$. 
(f) State $S3_{dl-d}$ with negative amplitude slope. Additional conditions are as in (e).
Other details are as in Fig. \ref{fig2_Gaussian_states}. 
Simulation parameters: $\lambda = 1$, $\omega = \pi$, $S = 1$, $N = 989$.
}
\label{fig6_HC_states}
\end{figure*}
\begin{figure*}
\centerfloat
\epsfig{figure=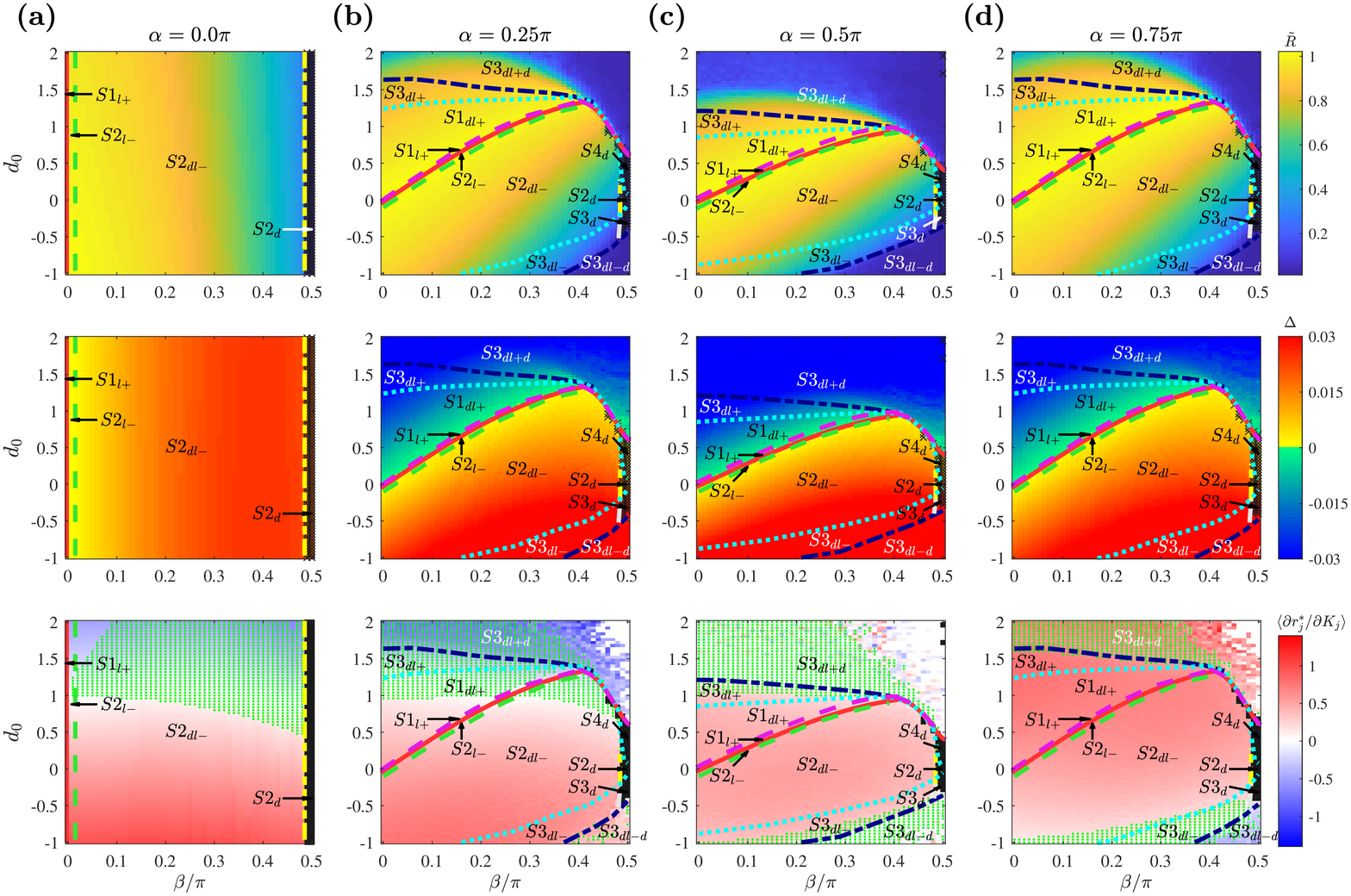, width=18.2cm}
\caption{Phase diagram with human brain network coupling strength distribution as a function of $\alpha$, $\beta$, and $d_0$ determining the form of coupling function. 
Representative fixed values of $\alpha$ are chosen at (a) $\alpha = 0$ where $\sin\alpha = 0,\, \cos\alpha = 1$, (b) $\alpha = 0.25\pi$ where $\sin\alpha > 0,\, \cos\alpha > 0$, (c) $\alpha = 0.5\pi$ where $\sin\alpha = 1,\, \cos\alpha = 0$, and (d) $\alpha = 0.75\pi$ where $\sin\alpha < 0,\, \cos\alpha > 0$.
Phase diagrams with order parameter $\tilde{R}$ (first row), phase pread $\Delta = \omega - \Omega$ (second row), and slope of the amplitude curve $(K_j, {r_j}^*)$ calculated as the average value of $\frac{\partial r^*_j}{\partial K_j}$ among the locked oscillators (third row). 
The boundaries and other details are as in Fig. \ref{fig3_Gaussian_diagram}.
For the obtained coupling strength set $K_j$ with $N = 989$, $(K_{\rm mean}, \sigma_K, K_{\rm min}, K_{\rm max}) = (36.5, 15.8, 1.01, 98.1) \times 10^{-3}$.
Simulation parameters: $\lambda = 1$, $\omega = \pi$, $S = 1$, $N = 989$.
}
\label{fig7_CIJ_diagram}
\end{figure*}
\begin{figure*}
\centering
\epsfig{figure=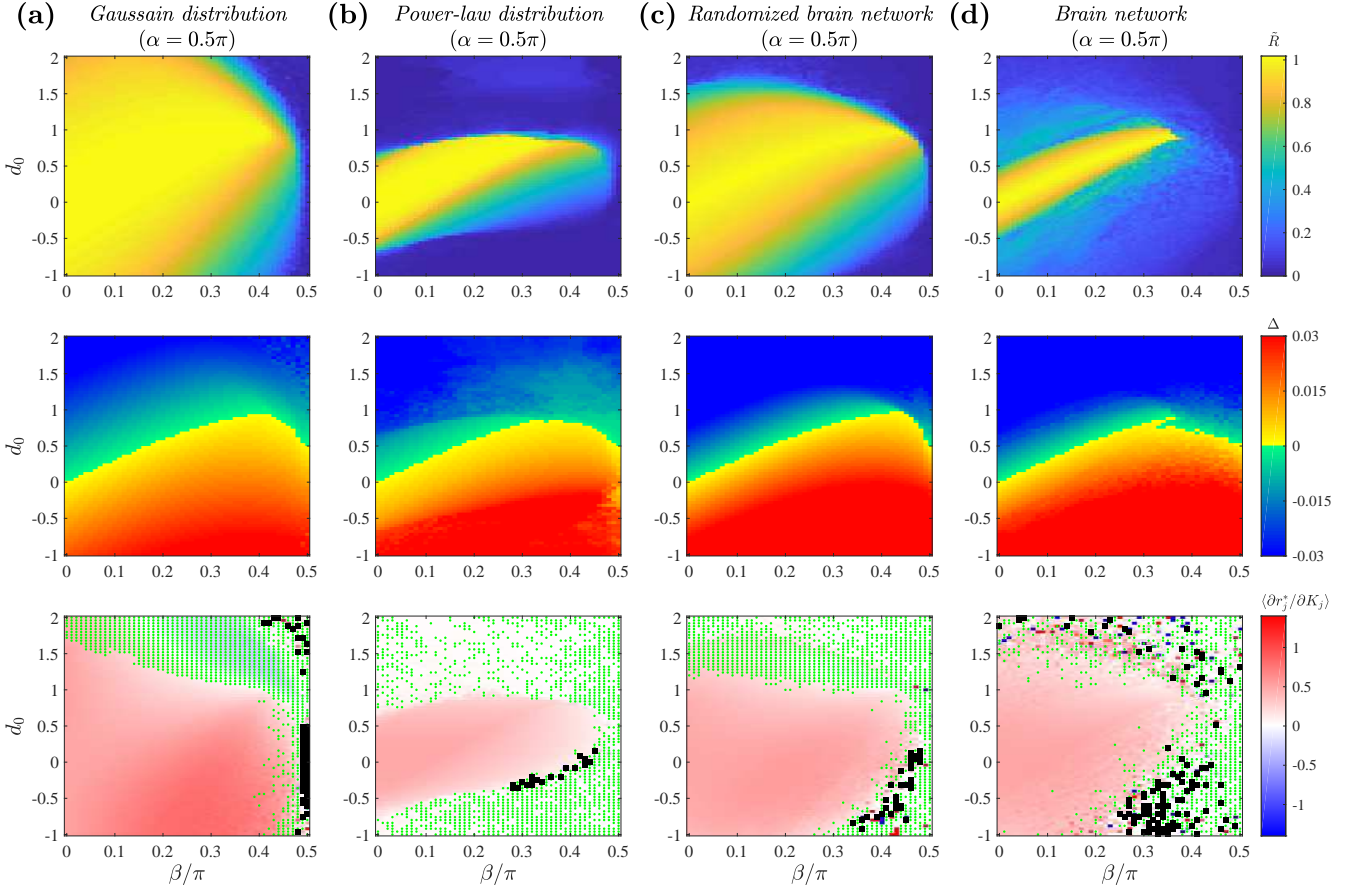, width=18.2cm}
\caption{Phase diagrams of complex networks using full connectivity matrices at $\alpha = 0.5\pi$. $\alpha$, $\beta$, and $d_0$ determining the form of coupling function. 
Phase diagrams at the representative fixed value of $\alpha = 0.5\pi$ are shown for (a) random network following Gaussian degree distribution, (b) random network following power-law degree distribution, (c) random network following brain network degree distribution, and (d) brain network with empirical connectivity profile.
Model dynamics are shown for order parameter $\tilde{R}$ (first row), $\Delta = \omega - \Omega$ (second row), and slope of the $(K_j, {r_j}^*)$ curve calculated as the average value of $\frac{\partial r^*_j}{\partial K_j}$ among the locked oscillators (third row).  In the third row, regions with fully drifting population are marked in black. The parameter space in which inflection point exists in the $(K_j, {r_j}^*)$ curve (from positive to negative slope) is additionally marked with green dots. 
Simulation parameters: $\lambda = 1$, $\omega = \pi$, $S = 1$, $N = 1000$ (Gaussian \& Power-law) or $N = 989$ (Brain networks). 
}
\label{fig8_full_diagram}
\end{figure*}
\begin{figure*}
\centering
\epsfig{figure=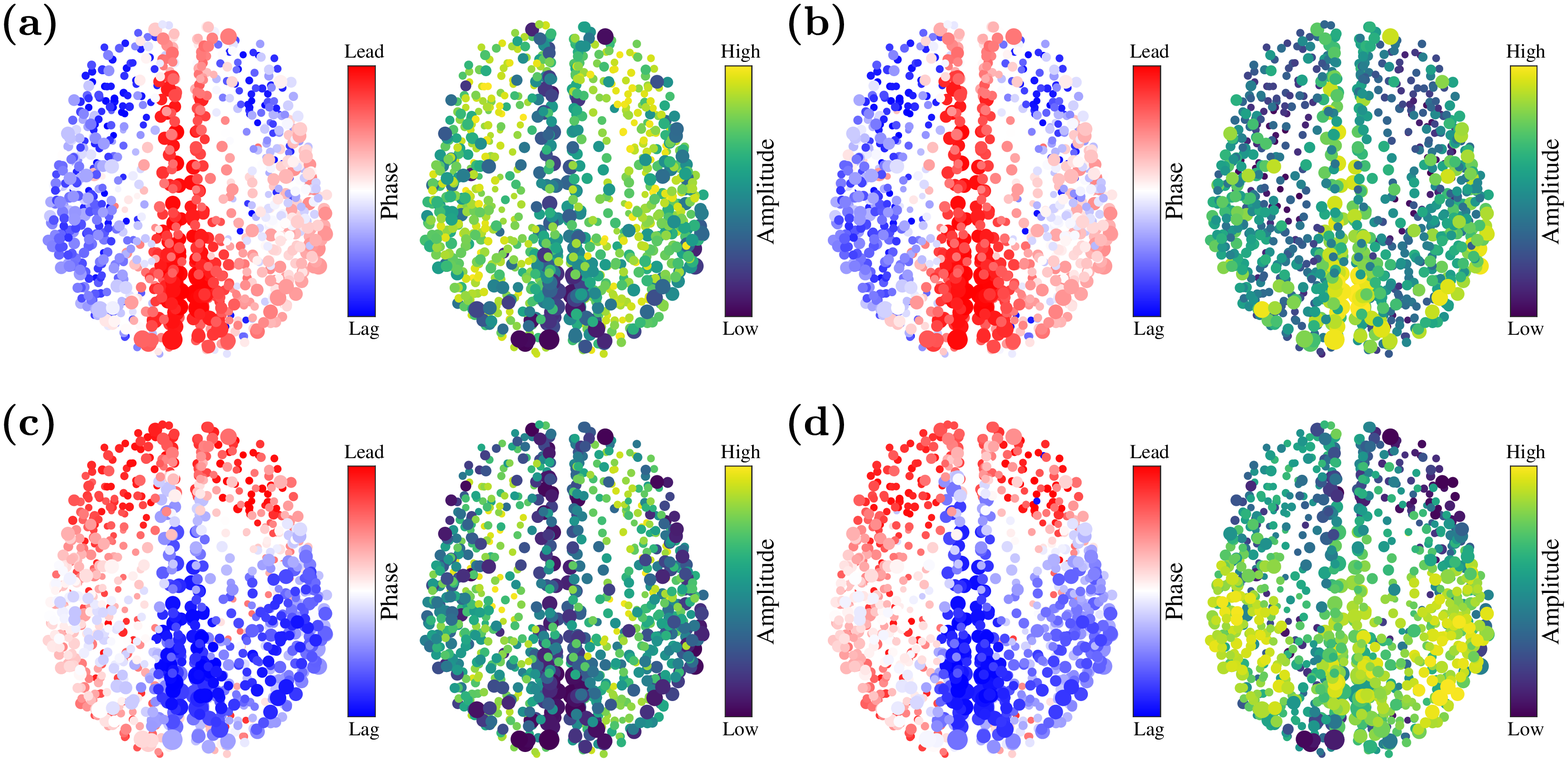, width=18cm}
\caption{Representative brain states, simulated from the model with full connectivity matrix (without mean-field approximation). (a) Brain states in which higher degree regions tend to phase-lead with lower amplitudes. ($\alpha$, $\beta$, $d_0$) $=$ ($0.25\pi, 0.22\pi,1.35$). (b) Higher degree regions phase-lead with higher amplitudes. ($\alpha$, $\beta$, $d_0$) $=$ ($0.25\pi, 0.21\pi,1.3$). (c) Higher degree regions phase-lag with lower amplitudes. ($\alpha$, $\beta$, $d_0$) $=$ ($0.1\pi, 0.2\pi,1.3$). (d) Higher degree regions phase-lag with higher amplitudes. ($\alpha$, $\beta$, $d_0$) $=$ ($0\pi, 0.1\pi,1.0$). 
Simulation parameters: $\lambda = 1$, $\omega = \pi$, $S = 1$, $N = 989$.
}
\label{fig9_brain_vis}
\end{figure*}

\section{\label{sec:simulations}Numerical Simulations}

We now describe numerical results of the oscillator model (Equation \ref{eq:SL1}) with coupling strengths derived from networks with different degree distributions: Gaussian, power-law, and brain-network-derived. All numerical simulations were carried out using a fourth order Runge-Kutta method with a fixed step size of $\Delta t = 0.01$. Oscillators have the identical intrinsic frequency $\omega=\pi$ ($0.5$ Hz) and $\lambda = 1$. For the initial conditions, each $\theta_j(0)$ was sampled randomly from $[0,2\pi)$ to form a near incoherent initial state. Each $r_j(0)$ was sampled randomly from a Gaussian distribution with mean $\sqrt{\lambda} = 1$ and standard deviation of  $0.1$. Unless noted otherwise, all results are averaged from ten different random initial conditions.

Fig. \ref{fig1_ComplexPlane} shows some representative examples of fully locked states (S1) on the complex plane after the stationary state is reached. These four states, which have different combinations of phase and amplitude dynamics, are all observable from model equation (\ref{eq:SL1}) in different parameter space. In Fig. \ref{fig1_ComplexPlane}(a), oscillators with larger coupling strengths (bigger and brighter dots) phase-lead oscillators with lower amplitude, whereas in (b) those with larger coupling strengths phase-lead with higher amplitude. In (c), oscillators with smaller coupling strengths (smaller and darker dots) phase-lead with lower amplitude, while (d) shows the opposite dynamics in amplitude. Drifting population is observable in other regions of parameter space, as shown in subsequent sections. 

\subsection{\label{subsec:random_dist}Gaussian coupling strength distributions}
{
First, we investigate the cases of Gaussian distribution for the coupling strength distributions with $N=1000$ oscillators. The values for inhomogeneous coupling strengths were randomly sampled from a Gaussian distribution with a mean of $20\times 10^{-3}$ and a standard deviation of $4.5\times 10^{-3}$. Figure \ref{fig2_Gaussian_states}(a) shows the distribution of the obtained coupling strength set $K_j$.

When coupling strengths were exactly or nearly homogeneous, or when $d_0\sin\alpha=\sin\beta$, the Stuart-Landau systems exhibited in-phase synchronous states. However, when coupling strengths were inhomogeneous and $d_0\sin\alpha \neq \sin\beta$, the system did not exhibit in-phase synchronous states.  

Representative examples of partially locked states of the system with the coupling strengths for various combinations of $\alpha$, $\beta$, and $d_0$ are presented in Figs. \ref{fig2_Gaussian_states}(b)-(f). 
%
Using the $\tilde{R}$ and $\Delta$ values obtained from the numerical simulation, we identified the self-consistent theoretical values of $(K_j,{\phi_j}^*)$ and $(K_j,{r_j}^*)$ from Eqs. (\ref{amp_eq}) and (\ref{eq:locked_phases1}). These are plotted against $K_j$ for the boundaries within the locking ranges (unshaded region in Figs. \ref{fig2_Gaussian_states}(b)-(f)). The analytically obtained values fit well with the simulations. As predicted by Eq. (\ref{eq:sign_slope_phi}), the signs of the slopes of the $(K_j,{\phi_j}^*)$ curve within the locking ranges are given by the negative of the sign of $\Delta$. Similarly, the signs of the slopes of $(K_j,{r_j}^*)$ curve are given by the sign of $\Delta$ in combination with ${\phi_j}^*$ and $\beta$. In Fig. \ref{fig2_Gaussian_states}(e) and (f), the inflection points in which the slope of the $(K_j,{r_j}^*)$ curve changes from positive to negative occurred when ${\phi_j}^*-\Phi+\beta = 0$ as predicted by Eq. (\ref{eq:sign_slope_phi}).

In order to understand the effects of the three parameters $\alpha$, $\beta$, and $d_0$ on the system, we plotted the phase diagrams by fixing $\alpha$ on a representative value, while varying $\beta$ and $d_0$. In this way, we tested different sign-combinations of the contribution from the constant term $d_0 e^{-i\alpha}$. For instance, at $\alpha = 0.5\pi$ the contribution of the constant term in the phase dynamics is maximal ($d_0\sin\alpha = d_0$ in Eq. 2), but it has no effect on the amplitude dynamics ($r_j d_0\cos\alpha = 0$ in Eq. 3). Conversely, at $\alpha = 0$, the $d_0 e^{-i\alpha}$ term affects only the amplitude dimension of the oscillators. Using this reasoning, four values of $\alpha$ at $\{0, 0.25\pi, 0.5\pi, 0.75\pi \}$ were selected. Setting $\alpha = 0.25\pi$ and $\alpha = 0.75\pi$ produces opposite-signed contributions in the $r_jd_0\cos\alpha$ term in Eq. (\ref{eq:SL3}).

Fig.\ref{fig3_Gaussian_diagram} shows the phase diagrams for the system with Gaussian coupling strength distributions. Each column depicts the phase diagram for a fixed $\alpha$ value, displayed with $\tilde{R}$, $\Delta$, and the average slope of the $(K_j, {r_j}^*)$ curve from the simulations. $\tilde{R}$ indicates the degree of global synchrony in the system, as used in Eq. (\ref{eq_sine_1}): the more locked oscillators in the system, the larger the $\tilde{R}$ values. $\Delta$ gives the information about the slope of the $(K_j, {\phi_j}^*)$ as shown in Eq. (\ref{eq:sign_slope_phi}). For example, larger magnitude of $\Delta$ indicate that the spread of phases among the locked oscillators is large, while the sign indicates whether high-degree or low-degree nodes are phase-leading. $\langle \frac{\partial r^*_j}{\partial K_j}\rangle$ gives the information about the slope of the $(K_j, {r_j}^*)$ curve. If $\langle \frac{\partial r^*_j}{\partial K_j}\rangle<0$, it suggests that higher-degree nodes are synchronized with lower amplitude at a stationary state. The boundaries between regions for states were determined from simulations of $\tilde{R}$ and $\Delta$ values. More details are provided in the caption of Fig. \ref{fig3_Gaussian_diagram}. 

The phase diagram for $\alpha = 0$ (Fig. \ref{fig3_Gaussian_diagram}.(a)) shows only four states are obtainable in the oscillator system: $S1_{l^+}$, $S2_{l^-}$, $S2_{dl^-}$, and $S2_{d}$. This is the fewest number of states observed for any value of $\alpha$. Note that the values of $\tilde{R}$ and $\Delta$ were nearly uniform across $d_0$ and were only modulated by the values of $\beta$, as indicated by Eq. (\ref{eq:SL2}). The vertical bounding curves also indicate that different states occurred only as a function of $\beta$ for $\alpha = 0$. Negative values of $\Delta$, for which state $S1_{l^+}$ occurs, were obtained near $\beta = 0$ (left of the red solid line). Further increases in the phase delay term $\beta$ produce increases in phase spread $\Delta$ and decreases in the order parameter $\tilde{R}$, until the system reaches an incoherent state ($S2_{d}$) near $\beta = \pi/2$ (right to the yellow dash-dotted line). The term, $d_0$, which describes the amplitude of the constant term applied to the system, only affects the slope of the ($K_j$, ${r_j}^*$) curve: for $d_0 > 1$, the system can exhibit synchronous states where, in the locked state, oscillators with larger coupling strengths have lower amplitudes.

With nonzero values of $\alpha$, additional sub-states in $S3$ and $S4$ were observed (Fig. \ref{fig3_Gaussian_diagram}.(b)-(d)). Of particular note, states in which oscillators with weaker coupling strength are locked while those with stronger coupling drift (e.g., $S3_{l^+d}$) were observed only at $\alpha \neq 0$. Such states are in contrast to our intuition that oscillators with larger coupling strength are easier to be locked.

The diagrams for $\tilde{R}$ and $\Delta$ (first and second row in Figure \ref{fig3_Gaussian_diagram}) reveal a symmetry around $\alpha = 0.5\pi$, such that $\alpha = 0.25\pi$ and $\alpha = 0.75\pi$ are nearly identical. This is as expected theoretically from the phase equation in Eq. (\ref{eq:SL2}), as the values of $\sin\alpha$ are symmetric around $\alpha = 0.5\pi$. This is not the case for $\cos\alpha$ in Eq. (\ref{eq:SL3}), and as expected, diagrams for the $(K_j, {r_j}^*)$ slope (third row) differ qualitatively. Yet the effects of the difference in the amplitude dynamics did not significantly influence the global synchrony of the system. This was so because coupling strengths $K_j$ normalized by $N$ result in a small overall coupling relative to the attraction to the limit cycle. This observation also explains the near-identity between the phase diagram for $\alpha = 0.5\pi$ (Fig. \ref{fig3_Gaussian_diagram}.(c)) with the one observed for the phase-reduced system of the model equation (\ref{eq:SL1}).\cite{chaos2019} At $\alpha = 0.5\pi$, the coupled amplitude term $r_j d_0\cos\alpha $ in Eq. (\ref{eq:SL3}) vanishes, and therefore the remaining amplitude dynamics were insufficient to affect the global synchrony of the system. (Note, however, that with the increase in the global coupling strength $S$, the difference in the amplitude dynamics is amplified; refer to Fig. S1 for more details.) 
}

\subsection{\label{subsec:power_law_dist}Power-law coupling strength distributions}
{
Next, we investigated a coupled Stuart-Landau system in which the coupling strengths were distributed according to a power-law distribution.  We keep the same number of oscillators ($N = 1000$). Fig. \ref{fig4_SF_states}(a) shows the coupling strengths distribution randomly sampled from a a truncated power-law distribution $P(x) \sim x^{-\gamma_0}$ with $\gamma_0 = 2$. Here, $K_{\rm min}$ and $K_{\rm max}$ were chosen such that the mean value of the distribution matches that of the Gaussian distributions used in the previous section.

Example observed states of the system are shown in Figures \ref{fig4_SF_states}(b)-(f). In contrast to the Gaussian coupling strength cases, in the power-law case there were only three types of partially locked states: $S2_{dl^-}$ (Fig. \ref{fig4_SF_states}(b)), $S3_{l^+d}$ (Fig. \ref{fig4_SF_states}(c)-(d)), and $S3_{dl^+d}$ (Fig. \ref{fig4_SF_states}(e)-(f)). Partially locked states $S1_{dl^+}$ and $S3_{dl^+}$  were not observed from the simulations for the power-law distribution, as seen in the phase diagrams in Fig. \ref{fig5_SF_diagram}. The most noticeable difference from the previous Gaussian case is the overall decrease in the synchronous region, indicated by the darker colored regions in the diagram for order parameter $\tilde{R}$ (Fig. \ref{fig5_SF_diagram}, top row). This suggests that in general, it is more difficult to achieve global synchrony for power-law distribution than for Gaussian distribution given the same parameter space. Although the two distributions have the same mean coupling strength $K_j$, the majority of oscillators in power-law distribution were assigned minimal coupling strength of $K_j < 0.01$. Therefore, the majority of oscillators contribute a smaller proportion of the total sum of coupling strengths in the power-law degree distribution, relative to the Gaussian degree distribution.

For nonzero cases of $\alpha$ (Fig. \ref{fig5_SF_diagram}.(b)-(d)), we also notice the absence of the regions for the partially locked states $S1_{dl^+}$ and $S3_{dl^+}$. However, a new partially locked state $S3_{dl^-d}$ was observed along the more negative values of $d_0$. In this state, only few oscillators with intermediate coupling strengths are locked. This is reflected in the lower values of $\tilde{R}$ for the $S3_{dl^-d}$ region.

The phase diagram at $\alpha = 0$ (Fig. \ref{fig5_SF_diagram}.(a)) shows the same pattern as with the previous section, where only $\beta$ affects the values of $\tilde{R}$ and $\Delta$. The major difference is that incoherent states are reached for lower values of $\beta$. However, the parameter region for fully drifting populations ($S2_d$) is smaller compared to the the Gaussian distributed couplings. This is reflected by the larger $K_{\rm max}$ value for the power-law coupling strength set, as the condition for $S2_d$ requires $K_{\rm max} \leq \frac{\Delta\cdot {r_{\rm max}}^*}{\tilde{R}}$ at $\alpha = 0$. Likewise, the smaller $K_{\rm min}$ value for the power-law coupling strength compared to the Gaussian case results in a decreased area for $S1_{l+}$ region, which requires $\frac{|\Delta|{r_{\rm min}}^*}{\tilde{R}} < K_{\rm min}$. Negative values of $\Delta$ were obtained only at $\beta = 0$. 
}

\subsection{\label{subsec:brain_net_dist}Brain network coupling strength distributions}
{
As an example model application to the real-world complex network, we now investigate oscillator systems in which the coupling strength approximates that of the couplings between regions of the human cerebral cortex. The distribution was derived from the network with 998 cortical regions \cite{honey2009}, where the coupling strength corresponds to the degree of each node in the network normalized by the total number of regions. Fig. \ref{fig6_HC_states}(a) shows the resulting coupling strength set $K_j$ with $N = 989$, after removing nodes with zero in-degree. 

Fig. \ref{fig7_CIJ_diagram} shows that the regions with fully locked population ($S1_{l+}$ and $S2_{l-}$) are significantly reduced, suggesting that fully locked states are practically non-existent in the system with brain network distribution. This phenomenon is accounted by the near-zero minimum coupling strength ($K_{\rm min} = 1.01\times 10^{-3}$) in the given network, compared to Gaussian and power-law distribution: with small values of $K_{\rm min}$, the parameter space with  $\frac{|\Delta|r^*_j}{\tilde{R}+D_0} < K_{\rm min}$ (for $S1_{l+}$) or with $\frac{\Delta r^*_j}{\tilde{R} \pm D_0} < K_{\rm min}$ (for $S2_{l-}$) decreases. Accordingly, the most commonly observed state in the phase diagrams was the partially locked state $S2_{dl-}$ (satisfying $K_{\rm min} \leq \frac{\Delta r^*_j}{\tilde{R} \pm D_0}$ with positive value of $\Delta$). 

Overall, we notice that the shape of phase diagrams for the brain network distribution is qualitatively in between the shapes derived for Gaussian and power-law distributions. The diagrams for $\tilde{R}$ (Fig. \ref{fig7_CIJ_diagram}, top row) show that the synchronous region is smaller than for the Gaussian case (Fig. \ref{fig3_Gaussian_diagram}) but larger than for the power-law case (Fig. \ref{fig5_SF_diagram}). The coupling strength distribution derived from large-scale brain networks (Fig. \ref{fig6_HC_states}(a)) also reflects the intermediate property between Gaussian and power-law distributions, as the majority of oscillators are skewed slightly left to $K_{\rm mean} = 0.0365$ and right to $K_{\rm min}$. (For Gaussian distribution, the majority of coupling strengths are centered around its $K_{\rm mean}$; for power-law distribution, around its $K_{\rm min}$).
}

\section{\label{sec:fullConnectivity}Simulation results with full complex networks connectivity}
{
The previous sections examined a mean-field model with the inhomogeneous coupling strength set $K_j$. In this section, we compare the simulation results of the coupled oscillator system on a full complex network. In the simulations, we retain the full connectivity profile, and so the simulated model is described by
\begin{eqnarray} \label{eq:full}
\dot{z}_{j} &=& \{\lambda_j - |z_{j}|^2 + i\omega_j \}z_{j} + \frac{S}{N}\sum_{k=1}^{N}A_{jk}(z_{k}e^{-i\beta}- z_{j}d_{0}e^{-i\alpha}), \nonumber \\
&&  j=1,2,...,N, ~\alpha \in [0,\pi) ,~\beta \in [0,\pi/2), ~d_0 \in \mathbb{R}, \nonumber \\
\end{eqnarray}
where $A_{jk}$ is the adjacency matrix describing the topology of the network. If $k$ influences $j$, we take $A_{jk} = 1$ and $A_{jk} = 0$ otherwise. With sufficiently large $N$, we can use the following mean-field approximation as in Ref. \onlinecite{twko2008partial}:
\begin{eqnarray} \label{eq:meanfield}
\sum_{k=1}^{N}A_{jk}H(z) \approx \frac{k_j}{N}\sum_{k=1}^{N}H(z),
\end{eqnarray}
where $k_j$ is the degree of oscillator $j$ and $H$ is the coupling function. Then, as an approximation of the model in Eq. (\ref{eq:full}), we can write 
\begin{eqnarray} \label{eq:SLapprox}
\dot{z}_{j} &=& \{\lambda_j - |z_{j}|^2 + i\omega_j \}z_{j} + \frac{S k_j}{N^2}\sum_{k=1}^{N}(z_{k}e^{-i\beta}- z_{j}d_{0}e^{-i\alpha}), \nonumber \\
\end{eqnarray}
which is equivalent to Eq. (\ref{eq:SL1}) with $K_j = \frac{k_j}{N}$. Thus the network characteristics are incorporated through the coupling inhomogeneity in $K_j$, a quantity directly proportional to the degree $k_j$. Eq. (\ref{eq:SL1}) is an approximation of Eq. (\ref{eq:full}) in the sense that the former treats the connections between nodes as all-to-all but normalize the effect to each node by its respective degree; hence it leaves out the topological information contained in $A_{jk}$. Yet such a mean-field method renders the system analytically tractable and amenable to numerical simulation, allowing for the kind of analysis that we have seen in the previous section. Thus, in this section we check whether the mean-field approximation qualitatively captures the dynamics obtained on the complex network topology. 

We simulate Eq. (\ref{eq:full}) with the complex networks investigated in the previous sections. Four networks were generated as follows, using the network generation algorithm provided in Ref. \onlinecite{bounova2015}. First, $N = 1000$ random positive integers $k_j$ were selected from the Gaussian distribution with mean $20$ and a standard deviation of $4.5$. Note that this is the same distribution used in the previous section (before deviding each $k_j$ by $N$). The randomly selected integers were bounded by $[k_{\rm min}=8,\, k_{\rm max}=34]$, which were determined by rounding $K_{\rm min}*N$ and $K_{\rm max}*N$ to the nearest integer, respectively. After the degree set was generated, each oscillator $j$ was randomly assigned $k_j$ neighbors without self-coupling, in such a way that the network is bidirectional. Second, using the same method, the power-law distribution with mean $20$ was used to generate a full network of size $N = 1000$ with $[k_{\rm min} = 6,\, k_{\rm max} = 121]$. Lastly for the brain network,
we use the provided brain network which retains the empirical connection data between brain regions. As with the previous section, nodes with zero in-degrees were omitted from the network. The resulting brain-derived graph contained  $N=989$ nodes. 

However, given the assumption of mean-field approximation which posits all-to-all connection (before normalization by degree), we expect that the unique connectivity profile of the empirical human brain network would not be incorporated. For this we reason we add another network following brain network degree distributions, but one in which the edges are randomized. This network was generated following the same method as above, where $N = 989$ random integers were selected from the degree distribution set bounded by $[k_{\rm min} = 1,\, k_{\rm max} = 97]$. 

For simplicity in presentation, we look at the phase diagrams for model equation (\ref{eq:SL1}) in the parameter space of $\alpha = 0.5\pi$ only. Fig. \ref{fig8_full_diagram} shows the simulation results for each given network. In the analogous figures from previous sections of this manuscript, we were able to mark the state boundaries; however, this was not possible for Figure because Eq. (\ref{eq:full}) does not allow for the ready classification of states as was possible for Eq. (\ref{pstate_conditions1})-(\ref{pstate_conditions2}).

Comparing the results in Fig. \ref{fig8_full_diagram}(a),(b) with those of the mean-field approximation in Fig.\ref{fig3_Gaussian_diagram}(c), \ref{fig5_SF_diagram}(c), respectively, we see that the phase diagrams show good agreement. For the case of brain-network-derived distribution seen in Fig \ref{fig7_CIJ_diagram}(c), we note that the results show better agreement with the brain network with randomized edges (Fig. \ref{fig8_full_diagram}(c)) than with the empirical brain network (Fig. \ref{fig8_full_diagram}(d)) as expected.

The empirical implication of Eq. (\ref{eq:SL1}) is that, in principle, various combinations of phase and amplitude dynamics are possible. Fig. \ref{fig9_brain_vis} shows example brain states visualized on the human cortical network, showing the same four combinations of increasing/decreasing amplitude and phase dynamics as in Fig. \ref{fig1_ComplexPlane}. For example, the state represented in Fig. \ref{fig9_brain_vis}(c) in which higher-degree regions phase-lag the rest of the network, while maintaining larger amplitudes, was empirically observed by Ref. \onlinecite{moon2015general}. In that work, it was shown that human subjects in eyes-closed resting states exhibit such behavior when time averaged over period of minutes in their EEG recording, as opposed to anesthetized unconscious state where the pattern disappears. Note that the parameter $\beta$ plays the role of time delays between neuronal regions; we estimate that the values will be around $0.1$ - $0.2$ in realistic brain simulations, in which each node is a cortical region.\onlinecite{moon2015general} In future, it may also be possible to empirically estimate the parameter $d_0$ by measuring changes in the peak-to-trough variation in neural oscillations.
}

\section{\label{sec:conclusion}Conclusions}
This study combined analytical and numerical methods to characterize the dynamics of coupled oscillators with both phase and amplitude dynamics. In particular, we set out to understand whether phase-locking occurred, and which nodes were leading and lagging in phase, depending on the form and strength of inter-node couplings. Using a mean-field model, we analytically mapped the effects of coupling strength inhomogeneity and the coupling functions, and documented the parameter regimes associated with phase-locked (synchronous and asynchronous), partially locked (e.g. partially drifting) and fully drifting states. The analytic results agreed with numerical simulations employing Gaussian distributions and power-law distributions for coupling strengths. In addition, we applied the model to understanding neural phase-lead and phase-lag relationships, by simulating and analyzing the model in the case where the node degrees are derived from empirical properties of cortical networks.  

These results will deepen the understanding of collective dynamics in complex systems.In particular, we determined conditions under which high-degree nodes can phase-lead or phase-lag the rest of the network, both when they have relatively higher amplitude and when they have lower amplitude. Furthermore, we showed that the high-degree nodes can have higher or lower amplitudes, regardless of whether they phase-lead or phase-lag the rest of the network. As a result, the system can exhibit four representative patterns as shown in Fig. \ref{fig9_brain_vis}. This finding is of particular significance in the modeling of neural systems, in which patterns of phase-leading and lagging along with the amplitude variations are associated with the control of information flow.\cite{palva2011localization,stam2012go,moon2015general}

Future studies could study the transient behavior of the model in the dynamics leading to the steady state, as the dynamics between stable states are important in complex real-world systems, such as brain dynamics. \cite{rabinovich2008transient, deco2012dynamical} For practical applications, it may also be useful to develop methods to identify, moment by moment, which sub-state within the phase-space (Figures \ref{fig3_Gaussian_diagram}, \ref{fig5_SF_diagram}, \ref{fig7_CIJ_diagram}) are occupied by real-world dynamical systems.


\section*{Acknowledgments}
We gratefully acknowledge the support of the National Institutes of Health (NIH), Bethesda, MD, USA (Grant No. RO1 MH111439).

\section*{Data Availability}
The data that support the findings of this study are available from the corresponding author upon reasonable request.

\section*{References}
\bibliographystyle{aipnum4-2}
\bibliography{Draft_revision}

\end{document}